\documentclass{article}
\usepackage[utf8]{inputenc}
\usepackage{amsmath}
\usepackage{amssymb}
\usepackage{amsthm}
\usepackage{graphicx}
\title{Solution of parameter-dependent diffusion equation in layered media}
\author{Antti Autio \and Antti Hannukainen\footnote{Corresponding author: e-mail: antti.hannukainen@aalto.fi} }
\date{\today}
\newtheorem{theorem}{Theorem}[section]
\newtheorem{lemma}{Lemma}[section]
\newtheorem{remark}{Remark}[section]

\newtheorem{assumption}{Assumption}[section]
\usepackage{subcaption}
\usepackage{color}
\usepackage{setspace}
\usepackage{tikz}
\usepackage{hyperref}

\begin{document}

\maketitle

\begin{abstract}
This work studies the parameter-dependent diffusion equation in a two-dimensional domain consisting of locally mirror symmetric layers. It is assumed that the diffusion coefficient is a constant in each layer. The goal is to find approximate parameter-to-solution maps that have a small number of terms. It is shown that in the case of two layers one can find a solution formula consisting of three terms with explicit dependencies on the diffusion coefficient. The formula is based on decomposing the solution into orthogonal parts related to both of the layers and the interface between them. This formula is then expanded to an approximate one for the multi-layer case. We give an analytical formula for square layers and use the finite element formulation for more general layers. The results are illustrated with numerical examples and have applications for reduced basis methods by analyzing the Kolmogorov $n$-width. 
\end{abstract}
\begin{keywords}
Parameter-dependent partial differential equations, low-rank approximation, Kolmogorov $n$-width
\end{keywords}\\\\
\begin{AMS}
35B30, 65N15, 65N30
\end{AMS}

\section{Introduction}
\begin{figure}
    \centering
    \begin{tikzpicture}
        \draw[black, very thick] (0,0) -- (2,0) -- (2,2) -- (1,1) -- (0,2) -- cycle
        node at (1,0.5) {$\Omega_1$};
        \draw[black, very thick] (2,0) -- (4,0) -- (4,2) -- (3,1) -- (2,2) -- cycle
        node at (3,0.5) {$\Omega_2$};
        \draw[black, very thick] (4,0) -- (6,0) -- (6,2) -- (5,1) -- (4,2) -- cycle
        node at (5,0.5) {$\Omega_3$};
    \end{tikzpicture}
    \caption{An example of a layered geometry with $N=3$. The diffusion coefficient is a constant $y_i$ in each subdomain $\Omega_i$.}
    \label{fig:geom_crown}
\end{figure}
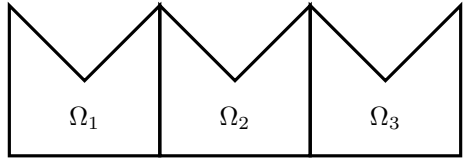

This work studies the diffusion equation in a layered domain $\Omega \subset \mathbb{R}^2$ consisting of $N$ locally mirror symmetric layers $\Omega_i$. We assume that the diffusion coefficient is a layer-wise constant with value $y_i$ in layer $\Omega_i$. The solution $u$ is viewed as a function from an admissible set of diffusion coefficient vectors $D \subset \mathbb{R}^N$ to the solution space $H^1_0(\Omega)$. Our goal is to approximate the parameter-to-solution map $u$ and study the approximation error. An example of the considered geometries is given in Fig. \ref{fig:geom_crown}. 

Our model problem is a simple example of parametric PDEs that arise e.g. in optimisation, solution of inverse problems, and uncertainty quantification, see \cite{CoDe:15,Su:15}. In recent years there has been interest in methods for computing approximate parameter-to-solution maps for such PDEs, see \cite{NoTeWe:08,QuMaNe:16}. The motivation for this work arises from authors' desire to understand the structure that in special cases allows some of these methods to perform exceptionally well, see \cite{KrTo:11,BaGr:15} and the discussion in \cite{MaPa:2021}. This phenomenon has been studied in \cite{BaCo:17} for $2\times 2$ chequerboard domains analytically and for $N\times N$ chequerboard domains numerically. Hence, only some cases with very simple parameter dependencies are currently well-understood.

Studying the special structure of the parameter-to-solution map requires one to develop tailored approximations to it that are not based on standard interpolation or approximation techniques. In this work, we develop such approximations from a novel solution formula for the two-layer case. For $N=2$, we show that 
\begin{equation} 
\label{eq:2x1formula}
u(y) = \frac{1}{y_1} w_{\Omega_1} + \frac{1}{y_2} w_{\Omega_2} + \frac{2}{y_1+y_2} w_\Gamma
\end{equation}
where functions $w_{\Omega_1},w_{\Omega_2},w_{\Gamma} \in H^1_0(\Omega)$ are independent of $y_1$ and $y_2$. 

Equation~\eqref{eq:2x1formula} is derived using a natural decomposition of $H^1_0(\Omega)$ into subdomain and interface spaces that are orthogonal in the energy inner product for any $y \in \mathbb{R}^2$, see \cite{BaCo:17}. Due to this orthogonality, the model problem splits to independent subdomain and interface problems. In \eqref{eq:2x1formula}, the intermediate  solutions related to the two subdomains are denoted as $w_{\Omega_1}$, $w_{\Omega_2}$ and the intermediate interface solution as $w_\Gamma$. The novelty of this work is an explicit solution formula for the interface problem. 

In the multi-layer case, i.e. $N>2$, we construct an approximate solution map \eqref{eq:Nformula} by utilizing the two-layer solution formula. We again use the natural subdomain-interface decomposition of the space $H^1_0(\Omega)$. In contrast to the two-layer case, the resulting interface problem is not analytically solvable. Instead, we further decompose the interface space into a finite dimensional \textit{slow} and an infinite dimensional \textit{fast subspace}. Then we prove that the interaction between functions in the fast subspaces related to different interfaces is small and thus the solution component from the fast subspace can be approximated by using the two-layer solution formula for each interface separately. The slow interface problem cannot be explicitly solved, but the solution is from a subspace $V_{s}$ whose dimension depends logarithmically on the desired accuracy of the approximation. 

The approximate solution for the multi-layer problem can then be expressed as
\begin{equation}
\label{eq:Nformula}
    u(y) \approx \sum_{i=1}^N \frac{1}{y_i} w_{\Omega_i} + \sum_{i=1}^{N-1} \frac{2}{y_i + y_{i+1}} w_{fi} + u_{s}(y)
\end{equation}
where $\{w_{\Omega_i} \}_{i=1}^N$ are the solutions to the parameter independent subdomain problems. The functions $\{w_{fi}\}_{i=1}^{N-1}$ are the solutions to parameter independent fast interface problems. The function $u_{s}(y)$ is the solution to the slow-interface problem that depends on the parameter vector $y$ and is posed in a subspace $V_{s}$ with small dimension. Therefore the given approximate parameter-to-solution map has $2N-1+\dim(V_{s})$ terms. The accuracy of this approximation as a function of the slow space dimension is analysed in detail.

The presented analysis is valid both in $H^1_0(\Omega)$ and a suitable finite element space. We explicitly give parameter-to-solution maps corresponding to rectangular layers. More general layers are treated by utilizing finite element discretisation of the spatial coordinate and using numerical techniques. The long-term goal of the authors is to utilize these techniques for efficient solution of the diffusion equation in an $N \times M$ chequerboard domain and to better understand the fundamental structure of this parametric equation. The rectangular $N\times 1$ case is a natural stepping stone towards that goal.

The approximate solution formula \eqref{eq:Nformula} is useful mostly for theoretical reasons.  It gives estimates for the Kolmogorov n-width and explains the success of sampling-based reduced basis methods for certain problems. In addition, it very concretely characterizes the structure of the solution. We have also included explanation on how \eqref{eq:Nformula} can be used in a computational setting for two reasons. Firstly, we believe that the numerical realisation helps to build intuition. Secondly, computational aspects are relevant for non-rectangular layers and for future extensions of this work to $N\times M$ chequerboard domains.

This article is organized as follows. First, we review the necessary background material. Then we derive the exact parameter-to-solution map in the two-layer case. Next, we treat the multi-layer case. We end this work with applications to reduced basis methods, numerical examples, and conclusions. 

\section{Model problem}

In this section, we first define the domain $\Omega$ and the layers $\{ \Omega_i \}_{i=1}^N$. Then we state the model problem, briefly discuss its properties, and introduce useful notation. 

Assume that the domain $\Omega \subset \mathbb{R}^2$ consists of $N$ layers $\Omega_i$. The interfaces between layers are vertical and the layers are locally mirror symmetric with respect to the interfaces. This is, each layer is obtained from layer $\Omega_1$ defined as 
\begin{equation*}
    \Omega_1 := \{ (x_1,x_2) \in \mathbb{R}^2 \; | \; x_1 \in (0,l) \; \mbox{and} \; x_2 \in (r_1(x_1),r_2(x_1)) \; \}
\end{equation*}
where $l>0$ is the layer thickness and $r_1, r_2 : (0,l) \mapsto \mathbb{R}$ are sufficiently regular functions satisfying $r_1(t) < r_2(t)$ for any $t\in (0,l)$. The remaining layers are then defined as 
\begin{equation}
\label{eq:mirrorsymm}
    \Omega_{i+1} = \{ \; (x_1,x_2)\in \mathbb{R}^2 \; | \;(2l i - x_1,x_2) \in \Omega_{i} \;\} \quad \mbox{for $i=1,\ldots,N-1$.}
\end{equation}
The interface between two layers is denoted by $\Gamma_i = \partial \Omega_{i} \cap \partial \Omega_{i+1}$ for $i=1,\ldots,N-1$ and the whole domain $\Omega$ is $\Omega = \text{int} \left( \cup_{i=1}^N \overline{\Omega}_{i} \right)$.

We now proceed to state the model problem. Let $f\in L^2(\Omega)$, $\mathbb{R}_+$ denote the set of non-negative real numbers, and $D \subset \mathbb{R}_+^N$ be the parameter set. Consider the problem: find the parameter-to-solution map $u: D \mapsto  H^1_0(\Omega)$ such that $u(y) \in H^1_0(\Omega)$ satisfies 
\begin{equation}
\label{eq:prob}
    \sum_{i=1}^N y_i (\nabla u(y),\nabla v)_{\Omega_i} = (f,v)_\Omega \quad \mbox{for any $v\in H^1_0(\Omega)$ and $y \in D$}. 
\end{equation}
Let $\omega \subseteq \Omega$. Here  $(\cdot,\cdot)_{\omega}$ denotes the $L^2(\omega)$ inner product of the arguments restricted to $\omega \subseteq \Omega$. Similarly the $L^2$-norm restricted to $\omega$ is denoted by $\| \cdot \|_\omega = (\cdot,\cdot)_\omega^{1/2}$. Assume that $D \subset (m,M)^N$ for some $m,M>0$. Then the problem \eqref{eq:prob} has a unique solution for any $y\in D$ by the Lax-Milgram lemma and it holds that 
\begin{equation*}
\| u(y)\|_{H^1_0(\Omega)} \leq m^{-1} \| f\|_\Omega \quad \mbox{for all $y \in D$}
\end{equation*}
where $\| u(y)\|_{H^1_0(\Omega)} := \| \nabla u(y) \|_\Omega$. Let the parameter dependent bilinear form $a: D \times H^1_0(\Omega) \times H^1_0(\Omega) \mapsto \mathbb{R}$ and linear functional $L:H^1_0(\Omega) \mapsto \mathbb{R}$ be defined as 
\begin{equation*}
    a(y;z,v) = \sum_{i=1}^N y_i (\nabla z,\nabla v)_{\Omega_i} \quad \mbox{and} \quad  L(v) = (f,v)_\Omega
\end{equation*}
for all $y \in D$ and $z,y \in H^1_0(\Omega)$. In the following, we do not explicitly emphasize the parameter dependency of the bilinear form, and write $a(z,v)$ instead of $a(y;z,v)$. Same convention is used for the parameter dependent energy norm $\| \cdot \|_E := a(y;\cdot,\cdot)^{1/2}$. 

For simplicity, we assume that $l=1$. This allows us to obtain much simpler formulas in the following sections without loss of generality.  

\section{Two-layer solution formula}

In this section, we derive an explicit formula for the parameter-to-solution map in the two-layer case. Recall that the domain $\Omega$ consists of two layers $\Omega_1$ and $\Omega_2$ that are mirror symmetric with respect to the vertical interface $\Gamma = \partial \Omega_1 \cap \partial\Omega_2$. First, define the spaces: 
\begin{equation}
\label{eq:def:V0i}
V_{\Omega_i} = \{ \; u \in H^1_0(\Omega) \; | \; \mathop{supp}(u) \subset \Omega_i \; \} \quad \mbox{for $i=1,2$}
\end{equation}
and
\begin{equation}
\label{eq:def:VG}
V_\Gamma := \{ \; u \in H^1_0(\Omega) \; | \; (\nabla u, \nabla v)_\Omega = 0 \quad \mbox{for any $v \in V_{\Omega_1} \cup V_{\Omega_2}$} \;\}.
\end{equation}
Intuitively speaking the spaces $V_{\Omega_1}$ and $V_{\Omega_2}$ are the zero extensions of $H_0^1(\Omega_1)$ and $H_0^1(\Omega_2)$ to $\Omega$ respectively. Any $v_\Gamma \in V_\Gamma$ is uniquely determined by its value on $\Gamma$ via the problem: Find $v_\Gamma \in H^1_0(\Omega)$ such that 
\begin{equation}
\label{eq:vgamma_harm_ext}
\begin{aligned}
-\Delta v_\Gamma &= 0 \quad \mbox{in $\Omega_1$ and $\Omega_2$} \\
v_\Gamma &= v_\Gamma|_\Gamma \quad \mbox{on $\Gamma$} \\
v_\Gamma &= 0 \quad \mbox{on $\partial \Omega$}.
\end{aligned}
\end{equation}

We can make two observations on $V_{\Omega_1}$, $V_{\Omega_2}$ and $V_\Gamma$. Firstly, these spaces are orthogonal both in the $(\nabla \cdot ,\nabla \cdot)_\Omega$-inner product and the parameter dependent inner product $a(y;\cdot,\cdot)$ for any $y\in D$. Secondly, it holds that $V_{\Omega_1}$, $V_{\Omega_2}$ and $V_\Gamma$ span $H^1_0$. To see this, let $v \in H^1_0$ and $v_\Gamma \in V_\Gamma$ have the same interface value on $\Gamma$, i.e., $v|_\Gamma = v_\Gamma|_\Gamma$. Then $(v-v_\Gamma)|_\Gamma = 0$ and thus $(v-v_\Gamma)|_{\Omega_i} \in H^1_0(\Omega_i)$ for $i\in \{1,2\}$. It follows that $H^1_0(\Omega) = V_{\Omega_1} \oplus V_{\Omega_2} \oplus V_\Gamma$. Since the spaces $V_{\Omega_1}$, $V_{\Omega_2}$, and $V_\Gamma$ are orthogonal and together span $H^1_0(\Omega)$, the problem \eqref{eq:prob} splits into three independent problems which are posed in $V_{\Omega_1}$, $V_{\Omega_2}$, and $V_\Gamma$ respectively. 

The problems in $V_{\Omega_1}$ and $V_{\Omega_2}$ amount to solving the Poisson equation in $\Omega_1$ and $\Omega_2$ with zero boundary conditions. The problem posed in $V_\Gamma$ is solved analytically in Theorem~\ref{thm:2x1} with the help of the following Lemma that relies on mirror symmetry of the domain $\Omega$ and the functions in $V_\Gamma$. Namely, by \eqref{eq:mirrorsymm} and \eqref{eq:vgamma_harm_ext} it holds that 
\begin{equation}
\label{eq:symm}
\begin{aligned}
v_\Gamma(1+t,x_2) &= v_\Gamma(1-t,x_2) \\
-\partial_x v_\Gamma(1+t,x_2) &= \partial_x v_\Gamma(1-t,x_2) \\
\partial_y v_\Gamma(1+t,x_2) &= \partial_y v_\Gamma(1-t,x_2)
\end{aligned}
\end{equation}
for any $v_\Gamma \in V_\Gamma$ and any $t \in (0,1), x_2$ for which $(1\pm t,x_2) \in \Omega$. 
\begin{lemma} 
\label{lemma:average} Let $V_\Gamma$ be as defined in \eqref{eq:def:VG}. Then there holds that
\begin{equation}
( y_1 \nabla u_\Gamma, \nabla v_\Gamma)_{\Omega_1} + (y_2 \nabla u_\Gamma, \nabla v_\Gamma)_{\Omega_2} = \frac{y_1+y_2}{2} (\nabla u_\Gamma, \nabla v_\Gamma)_{\Omega} 
\label{eq:gamma_split}
\end{equation} 
for any $u_\Gamma,v_\Gamma \in V_\Gamma$ and $(y_1,y_2) \in D$. 
\end{lemma}
\begin{proof}
Let $\overline{y}$ denote the mean of $y_1$ and $y_2$. Without loss of generality, assume that $y_1 > y_2$. Split the term $(\nabla u_\Gamma, \nabla v_\Gamma)_{\Omega} = (\nabla u_\Gamma, \nabla v_\Gamma)_{\Omega_1} + (\nabla u_\Gamma, \nabla v_\Gamma)_{\Omega_2}$. By simple manipulations, \eqref{eq:gamma_split} is equivalent to
\begin{equation*}
(y_1 - \overline{y}) ( \nabla u_\Gamma, \nabla v_\Gamma)_{\Omega_1} + (y_2 - \overline{y}) ( \nabla u_\Gamma, \nabla v_\Gamma)_{\Omega_2} = 0 \quad \mbox{for any $u_\Gamma,v_\Gamma \in V_\Gamma$}.
\end{equation*} 
Dividing by $(y_1 - \overline{y})$ gives 
\begin{equation*}
( \nabla u_\Gamma, \nabla v_\Gamma)_{\Omega_1} - ( \nabla u_\Gamma, \nabla v_\Gamma)_{\Omega_2} = 0 \quad \mbox{for any $u_\Gamma,v_\Gamma \in V_\Gamma$}. 
\end{equation*} 
The proof follows by observing that the above equation holds true by mirror symmetry of $u_\Gamma$ and $v_\Gamma$. Note that
\begin{equation*}
( \nabla u_\Gamma, \nabla v_\Gamma)_{\Omega_i} = \int_{\Omega_i} \partial_x u_\Gamma \partial_x v_\Gamma +  \int_{\Omega_i} \partial_y u_\Gamma \partial_y v_\Gamma \quad \mbox{for $i=1,2$}.
\end{equation*}
Application of \eqref{eq:symm} gives $( \nabla u_\Gamma, \nabla v_\Gamma)_{\Omega_2} = ( \nabla u_\Gamma, \nabla v_\Gamma)_{\Omega_1}$. 
\end{proof}

\begin{theorem} \label{thm:2x1} Let the spaces $V_{\Omega_1},V_{\Omega_2},V_\Gamma$ be as defined in \eqref{eq:def:V0i} and \eqref{eq:def:VG}. In addition, let $u:\mathbb{R}^2_+ \mapsto H^1_0(\Omega)$ be s.t. $u(y) \in H^1_0(\Omega)$ satisfies
\begin{equation*}
    y_1 (\nabla u(y),\nabla v)_{\Omega_1} + y_2 (\nabla u(y),\nabla v)_{\Omega_2} = (f,v)_\Omega \quad \mbox{for any $v\in H^1_0(\Omega)$}. 
\end{equation*}
Then 
\begin{equation*}
u(y) = \frac{1}{y_1} w_{\Omega_1} + \frac{1}{y_2} w_{\Omega_2} + \frac{2}{y_1+y_2} w_\Gamma,
\end{equation*}
where $w_{\Omega_1} \in V_{\Omega_1},w_{\Omega_2} \in V_{\Omega_2},w_\Gamma\in V_\Gamma$ satisfy
\begin{align}
\label{eq:N2_problem1}
    (\nabla w_{\Omega_1}, \nabla v_{\Omega_1})_\Omega & = (f,v_{\Omega_1})_\Omega \\
\label{eq:N2_problems2}
    (\nabla w_{\Omega_2}, \nabla v_{\Omega_2})_\Omega & = (f,v_{\Omega_2})_\Omega \\
    (\nabla w_\Gamma, \nabla v_\Gamma)_\Omega & = (f,v_\Gamma)_\Omega 
    \label{eq:N2_problems3}
\end{align}
for any $v_{\Omega_1} \in V_{\Omega_1},v_{\Omega_2} \in V_{\Omega_2},v_\Gamma \in V_\Gamma$. 
\end{theorem}
\begin{proof} As $H^1_0(\Omega) = V_{\Omega_1} \oplus V_{\Omega_2} \oplus V_\Gamma$, function $u(y) \in H^1_0(\Omega)$ admits the expansion $u(y) = u_{\Omega_1}(y) + u_{\Omega_2}(y) + u_{\Gamma}(y)$ for $u_{\Omega_1}(y) \in V_{\Omega_1}$, $u_{\Omega_2}(y) \in V_{\Omega_2}$, $u_\Gamma(y) \in V_\Gamma$. Using $a$-orthogonality of $V_{\Omega_1}$, $V_{\Omega_2}$ and $V_\Gamma$, these functions satisfy:
\begin{align*}
y_1 (\nabla u_{\Omega_1}(y) ,\nabla v_{\Omega_1} )_{\Omega} & = (f,v_{\Omega_1})_{\Omega} \\
y_2 (\nabla u_{\Omega_2}(y) ,\nabla v_{\Omega_2})_{\Omega} & = (f,v_{\Omega_2})_{\Omega} \\
y_1 (\nabla u_\Gamma(y) ,\nabla v_\Gamma)_{\Omega} + y_2 (\nabla u_\Gamma(y) ,\nabla v_\Gamma)_{\Omega} & = (f,v_\Gamma)_\Omega
\end{align*}
for any $v_{\Omega_1} \in V_{\Omega_1}$, $v_{\Omega_2} \in V_{\Omega_2}$ and $v_\Gamma \in V_{\Gamma}$. Application of Lemma~\ref{lemma:average} completes the proof. 
\end{proof}
\begin{remark} The space $V_\Gamma$ is non-standard. Thus, the easiest way to compute $w_\Gamma$ is not to solve~\eqref{eq:N2_problems3}, but to evaluate $u(y)$ at three appropriately chosen values of $y\in D$. The solution components can then be recovered by solving a $3\times 3$ linear system using \eqref{eq:2x1formula}.
\end{remark}

\subsection{Explicit solution of the interface problem}

In this section, we demonstrate how $w_\Gamma \in V_\Gamma$ defined in Equation~\eqref{eq:N2_problems3} is solved analytically for square layers. As problems \eqref{eq:N2_problem1} and \eqref{eq:N2_problems2} are standard Poisson equations in squares $\Omega_1,\Omega_2$ with zero Dirichlet boundary conditions, their solution is straightforward and not discussed. 

For simplicity, in this section we set $\Omega = (-1,1)\times(0,1)$, $\Omega_1 = (-1,0)\times(0,1)$, and $\Omega_2 = (0,1)^2$. Consider the following problem: for $n\in \mathbb{N}^+$ find $\psi_{n} \in V_\Gamma$ satisfying 
\begin{align*}
    -\Delta \psi_n &= 0, &&\mbox{in} \; \Omega_1 \cup \Omega_2 \\
    \psi_n &= \sin(\pi n x_2), &&\mbox{on} \; \Gamma \\
    \psi_n &= 0, &&\mbox{on} \; \partial \Omega
\end{align*}
in the weak sense. The solution to this problem is 
\begin{equation}
    \psi_{n} = \frac{1}{\sinh(\pi n)}\sin(\pi n x_2)\sinh(\pi n(1 - |x_1|))
    \label{eq:anal_harm_ext}
\end{equation}
for any $n\in \mathbb{N}^+$. 

We use this solution to characterise the space $V_ \Gamma$. Recall that any function in $L^2(\Gamma)$ can be expressed as a sine series and that the trace of any $H_0^1(\Omega)$ function on $\Gamma$ belongs to $H_{00}^{1/2}(\Gamma)$. Since it holds that $H_{00}^{1/2}(\Gamma) \subset L^2(\Gamma)$, the interface value of any $v_\Gamma \in H^1_0(\Omega)$ has the expansion 
\begin{equation*}
    v_\Gamma|_\Gamma = \sum_{n=1}^\infty \alpha_n \sin(\pi n x_2).
\end{equation*}

Using \eqref{eq:vgamma_harm_ext}, the space $V_\Gamma$ is then spanned by the harmonic extensions of admissible functions:
\begin{equation*}
    V_\Gamma = \left\{ \sum_{n=1}^\infty \alpha_n \psi_n \; \bigg{|} \; \left\| \sum_{n=1}^\infty \alpha_n \psi_n \right\|_{H^1_0(\Omega)} < \infty \right\}.
\end{equation*}
Thus, any $v_\Gamma \in V_\Gamma$ has the expansion
\begin{equation}
\label{eq:uGexpansion}
    v_\Gamma = \sum_{n=1}^\infty \alpha_n \psi_{n} = \sum_{n=1}^\infty \frac{\alpha_n}{\sinh(\pi n)} \sin(\pi n x_2) \sinh(\pi n(1 - |x_1|)).
\end{equation}

Next, we use the expansion \eqref{eq:uGexpansion} to solve $w_\Gamma$ from the interface space problem \eqref{eq:N2_problems3}. Choosing test functions $v_\Gamma = \psi_{m}$ yields
\begin{equation}
  \left(\nabla \sum_{n=1}^\infty \alpha_n \psi_{n}, \; \nabla \psi_{m} \right)_\Omega
    = \left( f, \; \psi_{m} \right)_\Omega.
    \label{eq:boundary_innerprod}
\end{equation}
for any $m \in \mathbb{N}^+$. The gradient can be calculated directly from \eqref{eq:anal_harm_ext}:
\begin{equation*}
    \nabla \psi_{n} = 
\frac{\pi n}{\sinh(\pi n)} 
\begin{bmatrix} 
\pm \sin( \pi n x_2) \cosh(\pi n (1 \pm x_1))\\
\cos( \pi n x_2) \sinh(\pi n (1 \pm x_1)) 
\end{bmatrix}
\end{equation*}
where the plus signs are used in $\Omega_1$ and the minus signs in $\Omega_2$. Due to the orthogonality of the basis, terms in \eqref{eq:boundary_innerprod} where $n \neq m$ are zero. We calculate the inner product and get
\begin{equation}
    \alpha_m = \frac{4 \sinh^2(\pi m)}{ \pi m \sinh(2\pi m)} \int_\Omega f\psi_m
    \label{eq:boundary_sol}
\end{equation}
which can be plugged into \eqref{eq:uGexpansion} to find $w_\Gamma$.

As an illustrative example, consider the case of $f=\sin(\pi x_2)$. In this case the integrals in \eqref{eq:boundary_sol} can be analytically calculated. The result is:
\begin{equation*}
    \alpha_1 = \frac{2 \sinh^2(\pi) \tanh\left(\frac{\pi}{2} \right)}{\pi^2 \sinh(2\pi)} 
\end{equation*}
and $\alpha_m=0$ if $m>1$. Therefore we get an explicit solution for the problem:
\begin{equation*}
    w_\Gamma = \alpha_1 \psi_1 = \frac{2 \sinh(\pi) \tanh\left(\frac{\pi}{2} \right)}{ \pi^2 \sinh(2\pi)} \sin(\pi x_2) \sinh(\pi (1 - |x_1|)).
\end{equation*}

\section{Multi-layer case}
In this section, we construct an approximate parameter-to-solution map for \eqref{eq:prob} in the case of $N>2$. Identical to the two-layer case, we begin with the subdomain-interface decomposition of $H^1_0(\Omega)$. For $i = 1,\ldots,N$ let $\Omega_{i}$ be as in \eqref{eq:mirrorsymm}, and define
\begin{equation*}
V_{\Omega_i} := \{ \; v \in H^1_0(\Omega) \; | \; \mathop{supp}(v) \subset \Omega_i \; \} 
\end{equation*}
as well as 
\begin{equation*}
V_{\Gamma} := \{ \; u \in H^1_0(\Omega) \; | \; (\nabla u, \nabla v)_\Omega = 0 \mbox{ for all $v \in V_{\Omega_j}$ and $j=1, \ldots, N$} \;\}.
\end{equation*}
Any $v_\Gamma \in V_\Gamma$ is uniquely determined by its value on the set of interfaces $\Gamma = \cup_{i=1}^{N-1} \Gamma_i$. Namely, each $v_\Gamma|_{\Omega_i}$ is the harmonic extension of $v_{\Gamma}|_{\Gamma_{i-1}}$ and $v_{\Gamma}|_{\Gamma_{i}}$ to $\Omega_i$, i.e., $v_\Gamma|_{\Omega_i} = w_i$, where $w_i \in H^1(\Omega_i)$ satisfies
\begin{equation*}
\left\{ 
\begin{aligned}
-\Delta w_i &= 0 &\mbox{in $\Omega_i$} \\
z_i &= v_\Gamma &\mbox{on $\Gamma_{i-1} \cup \Gamma_{i}$}
\end{aligned}
\right.
\end{equation*}
As harmonic extension is a linear operation, any $v_\Gamma \in V_\Gamma$ can be decomposed to a sum of local harmonic extensions $z_i \in V_\Gamma$ supported in $\Omega_i \cup \Omega_{i+1}$ and satisfying
\begin{equation}
\label{eq:V_Gamma_i_problem}
\left\{ 
\begin{aligned}
-\Delta z_i &= 0 &\mbox{in $\Omega_i$} \\
z_i &= 0 &\mbox{on $\partial \Omega_i \setminus \Gamma_{i}$} \\
z_i &= v_\Gamma &\mbox{on $\Gamma_{i}$}
\end{aligned}
\right.
\quad \quad 
\mbox{and} 
\quad \quad
\left\{ 
\begin{aligned}
-\Delta z_i &= 0 &\mbox{in $\Omega_{i+1}$} \\
z_i &= 0 &\mbox{on $\partial \Omega_{i+1} \setminus \Gamma_{i}$} \\
z_i &= v_\Gamma &\mbox{on $\Gamma_{i}$ }
\end{aligned}
\right. .
\end{equation}
Accordingly, the space $V_{\Gamma}$ admits the decomposition
\begin{equation*}
V_{\Gamma} = \bigoplus_{i=1}^{N-1} V_{\Gamma_i},
\end{equation*}
where
\begin{equation}
\label{eq:v_gamma_i}
V_{\Gamma_i} := \{ \; v \in V_\Gamma  \; | \;  v|_{\Gamma_j} = 0 \quad \mbox{for $j =1,\ldots,N-1,\; j\neq i$}  \; \} 
\end{equation}
are the spaces of local harmonic extensions. Each $z_i \in V_{\Gamma_i}$ is uniquely determined by its boundary value on $\Gamma_i$ by \eqref{eq:V_Gamma_i_problem}.

The spaces $\{V_{\Omega_i}\}^N_{i=1}$ and $V_\Gamma$ are orthogonal both in $H^1_0(\Omega)$ and energy inner products for any $y\in D$. Hence, the solution $u(y)$ to \eqref{eq:prob} satisfies 
\begin{equation}
    \label{eq:md:sub_int_deco}
    u(y) = \sum_{i=1}^N u_{\Omega_i}(y) + u_{\Gamma}(y), 
\end{equation}
where the \emph{interface solution} $u_\Gamma(y) \in V_{\Gamma}$ is such that
\begin{equation}
\label{eq:md:Gproblem}
a( u_{\Gamma}(y),v_{\Gamma}) = (f,v_{\Gamma}) \quad \mbox{for any $v_\Gamma \in V_\Gamma$} 
\end{equation}
and each $u_{\Omega_i}(y) \in V_{\Omega_i}$ is the solution to the subdomain problem
\begin{equation*}
a (u_{\Omega_i}(y),v_{\Omega_i}) = (f,v_{\Omega_i}) \quad \mbox{for any $v_{\Omega_i} \in V_{\Omega_i}$}.
\end{equation*}
For $i = 1,\ldots,N$ it clearly holds that
\begin{equation}
    \label{eq:md:sub_sol}
    u_{\Omega_i}(y) = y_i^{-1} w_{\Omega_i}
\end{equation}
where $w_{\Omega_i} \in V_{\Omega_i}$ satisfies 
\begin{equation*}
(\nabla w_{\Omega_i}, \nabla v_{\Omega_i}) = (f,v_{\Omega_i}) \quad \mbox{for all $v_{\Omega_i} \in V_{\Omega_i}$}. 
\end{equation*}
The topic of the rest of this section is approximate solution of \eqref{eq:md:Gproblem} by utilizing Lemma~\ref{lemma:average}. 

\subsection{Slow-Fast decomposition}
\label{sec:sfdeco}
We proceed to approximate the interface solution $u_\Gamma : D \mapsto V_\Gamma$. First $V_\Gamma$ is decomposed as a direct sum of subspaces that are approximately orthogonal in the $a(\cdot,\cdot)$-inner product. Then we solve \eqref{eq:md:Gproblem} separately in each of these subspaces, as if they would be exactly orthogonal in the energy-inner product. The approximate interface solution is obtained as a sum of these auxiliary solutions. The decomposition is chosen in such a way that most of the auxiliary problems can be solved using the two-domain solution formula.

In this section, we first define the decomposition of $V_\Gamma$. Then we give formulas for exact and approximate solution of \eqref{eq:md:Gproblem} related to this decomposition. Finally, we formally define what we mean by approximately orthogonal subspaces and relate the error between the exact and approximate solutions to the appropriate orthogonality-measure. 

Recall that $V_\Gamma$ has the decomposition
\begin{equation*}
V_{\Gamma} = \bigoplus_{i=1}^{N-1} V_{\Gamma_i},
\end{equation*}
where each $V_{\Gamma_i}$ is the space of local harmonic extensions as defined in \eqref{eq:v_gamma_i}. Functions in $V_{\Gamma_i}$ are defined by their value on $\Gamma_i$ via solving the problem \eqref{eq:V_Gamma_i_problem}. Hence, a decomposition of $V_{\Gamma_i}$ is obtained by decomposing the interface trace space $H^{1/2}_{00}(\Gamma_i)$ into two parts as $H^{1/2}_{00}(\Gamma_i) = Y_{si} \oplus Y_{fi}$. Then $V_{\Gamma_i} = V_{si} \oplus V_{fi}$, where
\begin{equation}
\begin{aligned}
\label{eq:sf_deco}
V_{fi} & := \{ \; v\in V_\Gamma  \; | \; v|_{\Gamma_i} \in Y_{fi} \mbox{ and } v|_{\Gamma_j} = 0\mbox{ for $j = 1,\ldots,N-1, j\neq i$ }  \; \} \\
V_{si} & := \{ \; v\in V_\Gamma  \; | \; v|_{\Gamma_i} \in Y_{si} \mbox{ and } v|_{\Gamma_j} = 0\mbox{ for $j = 1,\ldots,N-1, j\neq i$ }  \; \}.
\end{aligned}
\end{equation}
This is, $V_{fi}$ and $V_{si}$ are the spaces of local harmonic extensions of $Y_{fi}$ and $Y_{si}$, respectively. Our idea is to choose $Y_{si}$ and $Y_{fi}$ so that: 
\begin{itemize}
    \item[(i)] Corresponding spaces of local harmonic extensions, $V_{si}$ and $V_{fi}$, are $a(\cdot,\cdot)$-orthogonal for any $y\in D$.
    \item[(ii)] $V_{fi}$ is approximately $a(\cdot,\cdot)$-orthogonal to spaces of local harmonic extensions $V_{\Gamma_{i-1}}$ and $V_{\Gamma_{i+1}}$ related to adjacent interfaces for any $y\in D$.
\end{itemize}
In this section, we assume the existence of a decomposition satisfying (i) and (ii). Suitable decompositions are explicitly constructed in Sections~\ref{sec:sin} and \ref{sec:SVD}. Our intuitive idea is that the space $V_{fi}$ contains local harmonic extensions of rapidly oscillating functions on $\Gamma_i$ that decrease quickly when moving away from $\Gamma_i$. Due to this decrease, it is plausible that $V_{fi}$ satisfies (ii). Spaces $V_{si}$ contain slowly oscillating functions that do not have a similar decay property and thus they interact with slow functions related to other interfaces. We combine all of them into one slow subspace $V_s$ that is a direct sum of the slow interface subspaces 
\begin{equation}
\label{eq:slow_union}
V_{s} = \bigoplus_{i=1}^{N-1} V_{si}.
\end{equation}
We arrive to the \emph{slow-fast decomposition} 
\begin{equation}
\label{eq:deco}
    V_\Gamma = V_{s} \oplus V_{f1} \oplus \cdots \oplus V_{f (N-1)}. 
\end{equation}
The spaces $\{V_{fi}\}_{i=1}^{N-1}$ are called fast spaces and $V_{s}$ the slow space.

We proceed to explain how decomposition \eqref{eq:deco} is used to approximately solve \eqref{eq:md:Gproblem}. Decompose the exact interface solution as $u_\Gamma(y) = \sum_{i=1}^{N-1} u_{fi}(y) + u_{s}(y)$ for $u_{fi}(y)\in V_{fi}$ and $u_s(y) \in V_{s}$. By \eqref{eq:md:Gproblem}, each $u_{fi}(y)$ and $u_s(y)$ satisfy 
\begin{equation}
\label{eq:uexact}
\begin{aligned}
a(u_{fi}(y),v_{fi}) &= L(v_{fi}) - a(u_\Gamma(y)-u_{fi}(y), v_{fi}) \\ 
a(u_{s}(y),v_{s}) &= L(v_{s}) - a(u_\Gamma(y)-u_{s}(y), v_{s})    
\end{aligned}
\end{equation}
for any $v_{fi} \in V_{fi}$ and $v_s\in V_s$. Under assumption (i), the \emph{interaction terms} in \eqref{eq:uexact} satisfy
\begin{equation*}
\begin{aligned}
a(u_\Gamma(y)-u_{fi}(y), v_{fi}) & = a(u_{i-1},v_{fi}) + a(u_{i+1},v_{fi}) \\
a(u_\Gamma(y)-u_{s}(y), v_{s}) 
& = \sum_{i=1}^{N-1} \left[ a(u_{f(i-1)},v_{si}) + a(u_{f(i+1)},v_{si}) \right].
\end{aligned}
\end{equation*}
Here $u_{f0} = u_{fN} = 0$ and $u_{0} = u_N = 0$. If assumption (ii) holds, these terms are small. The approximate interface solution is obtained by neglecting them from \eqref{eq:uexact}. Let the \emph{approximate interface solution} related to the decomposition \eqref{eq:deco} be  
\begin{equation}
\label{eq:IF_error}
\tilde{u}_\Gamma(y) = \sum_{i=1}^{N-1} \tilde{u}_{fi}(y) + \tilde{u}_{s}(y)
\end{equation}
where $\tilde{u}_{fi}(y) \in V_{fi}$ satisfy
\begin{align}
\label{eq:tildeuf}
a(\tilde{u}_{fi}(y),v_{fi}) &= L(v_{fi}) \\ 
\label{eq:tildeus}
a(\tilde{u}_{s}(y),v_{s}) &= L( v_{s}) 
\end{align}
for any $v_{fi} \in V_{fi} \; \; i=1,\ldots,N-1$, $v_{s} \in  V_{s}$, and $y\in D$. 

We proceed to solve \eqref{eq:tildeuf} by using the two-domain solution formula. By definition~\eqref{eq:sf_deco} it holds that $\mathop{supp} v_{fi} \subset \Omega_i \cup \Gamma_i \cup \Omega_{i+1}$ for all $v_{fi} \in V_{fi}$. Thus, \eqref{eq:tildeuf} reads 
\begin{equation}
\label{eq:abu_fast}
y_i (\tilde{u}_{fi}(y),v_{fi})_{\Omega_i}
+ y_{i+1} (\tilde{u}_{fi}(y),v_{fi})_{\Omega_{i+1}}
= L(v_{fi})
\end{equation}
Application of Lemma~\ref{lemma:average} gives the solution 
\begin{equation*}
    \tilde{u}_{fi}(y) = \frac{2}{y_{i}+y_{i+1}} w_{fi},
\end{equation*}
where $w_{fi} \in V_{fi}$ is the solution to \eqref{eq:tildeuf} with $y=1$. The slow problem in \eqref{eq:tildeus} has to be solved for each parameter, but it is posed in a parameter independent space that has a small dimension and an explicitly known basis.
\begin{remark} A central assumption in this work is that $V_{si}$ and $V_{fi}$ are $a$-orthogonal and that each problem posed in $V_{fi}$ is solved using the two-domain solution formula. This requires that the domains are mirror symmetric. Extension of our work to non-symmetric case requires a different approach for solving \eqref{eq:abu_fast} as well as constructing spaces $V_{si}$ and $V_{fi}$. For example, one can resort to interpolation with respect to variables $(y_i,y_{i+1})$ when  solving~\eqref{eq:abu_fast} leading to slightly more complicated expression for $\tilde{u}_{fi}(y)$. The authors are working on modifying the analysis for approximately $a$-orthogonal spaces $V_{si}$ and $V_{fi}$ in the non-symmetric case. 
\end{remark}
\subsection{Error estimates}
\label{sec:error_estimates}

We proceed to estimate the error between exact and approximate interface solutions, $u_\Gamma(y)$ and $\tilde{u}_\Gamma(y)$. The given results are used in Section~\ref{sec:sin} to evaluate the error for slow-fast decomposition based on harmonic extension of sinusoidal interface functions for square layers and in Section~\ref{sec:SVD} for SVD based slow-fast decomposition related to finite element discretisation of \eqref{eq:prob}. 

The error estimates are inspired by literature on the Schwarz method, e.g.~\cite{ToWi:2004, SmBjGr:2004} and are based on perturbation analysis. Recall that by \eqref{eq:uexact} and \eqref{eq:tildeuf}, the approximate interface solution is obtained by neglecting the interaction terms $a(u_\Gamma(y)-u_{fi}(y), v_{fi})$, $a(u_\Gamma(y)-u_{s}(y), v_{s})$ from \eqref{eq:uexact}. Error estimates are derived by using a simple perturbation argument and a suitable estimate for the neglected terms. In addition, we must establish that the slow-fast decomposition is stable, see \cite{Os:1994}. 

Make a standing assumption that $V_{si}$ and $V_{fi}$ are $a$-orthogonal for any $y\in D$. In Section \ref{sec:sin}, we present an explicit construction of these spaces for rectangular layers in the continuous space and in Section \ref{sec:fem} for arbitrary mirror symmetric layers in the finite element space. Both of these constructions satisfy the criteria (i) and (ii) given in Section~\ref{sec:sfdeco}. 

Recall that any $v_{\Gamma_i} \in V_{\Gamma_i}$ satisfies $\mathop{supp}(v_{\Gamma_i}) \subset \Omega_i \cup \Gamma_i \cup \Omega_{i+1}$. We may now estimate the term $a(u_\Gamma-u_{fi}, v_{fi})$ by using the \emph{fast interaction} term 
\begin{equation}
\label{eq:f2e}
    \epsilon_f = \max_{i=1,\ldots,N-2} \max_{\substack{z_{fi} \in V_{fi} \\ z_{i+1} \in V_{\Gamma_{i+1} }}} \frac{ (\nabla z_{fi},\nabla z_{i+1})_{\Omega_{i+1}} }{ \|\nabla z_{fi} \|_{\Omega_{i+1}} \|\nabla z_{i+1} \|_{\Omega_{i+1}}}.
\end{equation}
The fast interaction term describes how much the fast subspace related to the interface $i$ interacts with the harmonic extensions from the neighbouring interface. This is an interaction which is disregarded in approximations~\eqref{eq:tildeuf} and~\eqref{eq:tildeus}. The size of the fast interaction term depends on the choice of $V_{fi}$ and can be made arbitrarily small.

Stability of the slow-fast decomposition is linked to the \emph{edge-to-edge interaction term}
\begin{equation}
\label{eq:e2e}
    \Lambda = \max_{i=1,\ldots,N-2} \max_{\substack{z_{i} \in V_{\Gamma_i} \\ z_{i+1} \in V_{\Gamma_{i+1}} }} \frac{ (\nabla z_i,\nabla z_{i+1})_{\Omega_{i+1}} }{ \|\nabla z_{i} \|_{\Omega_{i+1}} \|\nabla z_{i+1} \|_{\Omega_{i+1}}}.
\end{equation}
Here, $\Lambda$ is a constant that measures how much interaction there is between $V_\Gamma$-functions related to neighbouring interfaces. Clearly, it holds that $\Lambda < 1$, as $\Lambda = 1$ can only be attained if 
\begin{equation*}
(\nabla z_i,\nabla z_{i+1} )_{\Omega_{i+1}} = \| \nabla z_i \|_{\Omega_{i+1}} \|\nabla z_{i+1} \|_{\Omega_{i+1}},
\end{equation*}
that only holds for $z_i=z_{i+1}$, which is never true. 

We proceed to prove stability of the slow-fast decomposition. First, we establish stability of the decomposition $V_\Gamma = V_{\Gamma_1} \oplus \cdots \oplus V_{\Gamma_{N-1}}$ that we call \emph{$V_\Gamma$-stability}. 
\begin{lemma}[$V_\Gamma$-stability]  
\label{lemma:Gstab} For $i = 1,\ldots, N-1$ let $V_{\Gamma_i}$ be as defined in~\eqref{eq:v_gamma_i} and $w_{\Gamma_i} \in V_{\Gamma_i}$. Then there holds that
\begin{equation*}
(1 - \Lambda) \sum_{i=1}^{N-1} \left\| w_{\Gamma_i} \right\|^2_E \leq \left\| \sum_{i=1}^{N-1} w_{\Gamma_i} \right\|^2_E \leq (1 + \Lambda) \sum_{i=1}^{N-1} \| w_{\Gamma_i} \|^2_E 
\end{equation*}
\noindent for any $y \in D$. 
\end{lemma}
\begin{proof}
Using the definition of the energy norm and recalling that $\mathop{supp}{w_{i}} \subset \Omega_{i} \cup \Gamma_{i} \cup \Omega_{i+1}$ yields
\begin{equation}
\label{eq:proof:Lstab1}
\left\| \sum_{k=1}^{N-1} w_{\Gamma_i} \right\|_E^2   
=
\sum_{k=1}^{N-1} \| w_{\Gamma_i} \|_E^2   
+ 
2 \sum_{i=1}^{N-2} a(w_{\Gamma_i},w_{\Gamma_{i+1}})
\end{equation}
By definition of $V_{\Gamma_i}$, $V_{\Gamma_{i+1}}$, and the bilinear form $a$, 
\begin{equation*}
2|a(w_{\Gamma_i},w_{\Gamma_{i+1}})|  = 2 y_{i+1} \left| \frac{(\nabla w_{\Gamma_i}, \nabla w_{\Gamma_{i+1}})_{\Omega_{i+1}}}{\|\nabla w_{\Gamma_i} \|_{\Omega_{i+1}} \| \nabla w_{\Gamma_{i+1}} \|_{\Omega_{i+1}}} \right |\|\nabla w_{\Gamma_i} \|_{\Omega_{i+1}} \| \nabla w_{\Gamma_{i+1}} \|_{\Omega_{i+1}} 
\end{equation*}
Using the definition of the edge-to-edge interaction term in \eqref{eq:e2e} gives 
\begin{equation*}
2|a(w_{\Gamma_i},w_{\Gamma_{i+1}})|   \leq 2 y_{i+1} \Lambda \|\nabla w_{\Gamma_i} \|_{\Omega_{i+1}} \| \nabla w_{\Gamma_{i+1}} \|_{\Omega_{i+1}}.
\end{equation*}

By inequality $2ab \leq a^2 + b^2$, it holds that
\begin{equation*}
2|a(w_{\Gamma_i},w_{\Gamma_{i+1}})| \leq \Lambda \| y_{i+1}^{1/2} \nabla w_{\Gamma_i} \|^2_{\Omega_{i+1}} + \Lambda \| y_{i+1}^{1/2} \nabla w_{\Gamma_{i+1}} \|^2_{\Omega_{i+1}}.
\end{equation*}
So that
\begin{equation}
\label{eq:proof:Lstab2}
\left|2 \sum_{i=1}^{N-2} a(w_{\Gamma_i},w_{\Gamma_{i+1}}) \right| \leq \Lambda \sum_{i=1}^{N-1} \| w_{\Gamma_i} \|_E^2.
\end{equation}
%
Combining \eqref{eq:proof:Lstab1} and \eqref{eq:proof:Lstab2} completes the proof. 
\end{proof}
\begin{lemma}[Slow-fast stability] 
\label{lemma:stab}
Let $V_{s}$, $V_{si}$, and $V_{fi}$ be as defined in~\eqref{eq:sf_deco} and \eqref{eq:slow_union}. Assume that $V_{si}$ and $V_{fi}$ are $a$-orthogonal for $i = 1,\ldots, N-1$ and any $y \in D$. In addition, let $w_{s} \in V_{s}$ and $w_{fi} \in V_{fi}$ for $i = 1,\ldots,N-1$. Assume that $\epsilon_f$ is sufficiently small. Then there holds that
\begin{equation*}
\gamma^2 \left\| w_{s} + \sum_{i=1}^{N-1} w_{fi} \right\|_E^2 \geq \left\| w_{s} \right\|_E^2 + \sum_{i=1}^{N-1} \left\| w_{fi} \right\|_E^2,
\end{equation*} 
where $\gamma^2 = \max\{ 1-2 \epsilon_f(1-\Lambda)^{-1}, 1 -\Lambda-\epsilon_f\}^{-1}$. 
\end{lemma}

\begin{proof}

By definition of the energy norm, 
\begin{equation*}
\left\| w_{s} + \sum_{i=1}^{N-1} w_{fi} \right\|_E^2 = \| w_{s} \|_E^2 + \left\| \sum_{i=1}^{N-1} w_{fi} \right\|_E^2
+ 
\sum_{i=1}^{N-1} 2 a (w_{fi}, w_{s} ).
\end{equation*}
By~\eqref{eq:slow_union} $w_{s} = \sum_{i=1}^{N-1} w_{si}$ for $w_{si} \in V_{si}$. Using this expansion, $a$-orthogonality between $V_{si}$ and $V_{fi}$, and recalling that $\mathop{supp}w_{fi} = \mathop{supp} w_{si} = \Omega_{i} \cup \Gamma_i \cup \Omega_{i+1}$ yields
\begin{equation*}
\sum_{i=1}^{N-1} a (w_{fi}, w_s ) = \sum_{i=2}^{N-1} a( w_{fi}, w_{s(i-1)} )+ \sum_{i=1}^{N-2} a(w_{fi}, w_{s(i+1)} ).
\end{equation*}
By definition of the fast interaction term in \eqref{eq:f2e} and inequality $2ab \leq a^2 + b^2$,
\begin{align*}
2\left|a( w_{fi}, w_{s(i-1)} )\right| & \leq \epsilon_f \left\| y_i^{1/2} \nabla w_{fi} \right\|^2_{\Omega_{i}} + \epsilon_f \left\| y_i^{1/2} \nabla w_{s(i-1)} \right\|^2_{\Omega_{i}} \\
2\left|a( w_{fi}, w_{s(i+1)} )\right| & \leq \epsilon_f \left\| y_{i+1}^{1/2} \nabla w_{fi} \right\|^2_{\Omega_{i+1}} + \epsilon_f \left\| y_{i+1}^{1/2} \nabla w_{s(i+1)} \right\|^2_{\Omega_{i+1}}
\end{align*}
Thus,
\begin{equation*}
2 \left|\sum_{i=1}^{N-1} a (w_{fi}, w_s )\right|  \leq  2 \epsilon_f \sum_{i=1}^{N-1} \left\| w_{si} \right\|_E^2 + \epsilon_f \sum_{i=1}^{N-1} \| w_{fi} \|_E^2.
\end{equation*}
Application of $V_\Gamma$-stability result in Lemma~\ref{lemma:Gstab} gives
\begin{equation*}
2 \left|\sum_{i=1}^{N-1} a (w_{fi}, w_s )\right|  \leq  \frac{2 \epsilon_f}{1-\Lambda} \| w_s \|_E^2 + \epsilon_f \sum_{i=1}^{N-1} \| w_{fi} \|_E^2.
\end{equation*}
and
\begin{equation*}
\left\| \sum_{i=1}^{N-1} w_{fi} \right\|_E^2  \geq  (1-\Lambda) \sum_{i=1}^{N-1} \| w_{fi} \|_E^2.
\end{equation*}
So that
\begin{equation*}
\left\| w_s + \sum_{i=1}^{N-1} w_{fi} \right\|_E^2 \geq \left(1-\frac{2 \epsilon_f}{1-\Lambda}\right) \| w_s \|_E^2 + (1 - \Lambda - \epsilon_f) \sum_{i=1}^{N-1} \| w_{fi} \|_E^2
\end{equation*}
\end{proof}
\begin{remark} In the proof of Lemma~\ref{lemma:stab}, the $V_\Gamma$-stability result of Lemma~\ref{lemma:Gstab} is applied to bound the terms $\sum_{i=1}^{N-1} \| w_{si} \|_E^2$ and $\sum_{i=1}^{N-1} \| w_{fi} \|_E^2$. The bound related to $w_{si}$ cannot be improved whereas the bound related to $w_{fi}$ could be improved by using \eqref{eq:f2e}. Such a modification is omitted to keep the presentation technically simpler. 
\end{remark} 
We proceed to estimate error in the approximate slow and fast solutions, i.e. $\| u_{s} - \tilde{u}_{s} \|_E$ and $\|u_{fi} - \tilde{u}_{fi}\|_E$. 

\begin{lemma}[Slow error estimate]
\label{lemma:serror}
Let $V_{s}$, $V_{si}$, and $V_{fi}$ be as defined in~\eqref{eq:sf_deco} and \eqref{eq:slow_union}. In addition, let $\epsilon_f$ be as defined in \eqref{eq:f2e}. Assume that $V_{si}$ and $V_{fi}$ are $a$-orthogonal for $i = 1,\ldots, N-1$ and any $y\in D$.  Assume that $\epsilon_f$ is sufficiently small.  Let $u_s(y),\tilde{u}_s(y) \in V_s$ satisfy \eqref{eq:uexact} and \eqref{eq:tildeus}, respectively. Then there exists $C$ independent on $y$, $u$, and $\epsilon_f$ such that
\begin{equation*}
\| u_{s}(y) - \tilde{u}_{s}(y) \|_E 
\leq
C \epsilon_f \gamma \| u_\Gamma(y) \|_E \quad \mbox{for any $y\in D$}
\end{equation*}
where $\gamma$ is as defined in Lemma~\ref{lemma:stab} and $u_\Gamma(y)$ is the interface solution defined in \eqref{eq:md:Gproblem}.
\end{lemma}
\begin{proof} By \eqref{eq:md:Gproblem} and \eqref{eq:uexact}, it holds that 
\begin{equation*}
\begin{aligned}
a(u_{s}(y),v_{s}) &= L(v_{s}) - a(u_\Gamma(y)-u_{s}(y), v_{s})  \\
a(\tilde{u}_{s}(y),v_{s}) &= L(v_{s}).
\end{aligned}
\end{equation*}
Hence, $\| u_{s} - \tilde{u}_{s} \|^2_E 
= 
-a(u_\Gamma(y)-u_{s}(y), u_{s} - \tilde{u}_{s})$. Recalling that $u_\Gamma(y)-u_{s}(y) = \sum_{i=1}^{N-1} u_{fi}$ yields 
\begin{equation*}
\| u_{s} - \tilde{u}_{s} \|^2_E 
= \sum_{i=1}^{N-1} -a( u_{s} - \tilde{u}_{s}, u_{fi} ).
\end{equation*}
Let $e_s = u_{s} - \tilde{u}_{s}$ and split $e_{s} = e_{s1} + e_{s2} + \cdots + e_{s(N-1)}$ for $e_{si} \in V_{si}$. Next, we utilize the $a$-orthogonality between $V_{si}$ and $V_{fi}$ to estimate the cross term in the above equation. It holds that
\begin{equation*}
\| u_{s} - \tilde{u}_{s} \|^2_E 
=
\sum_{i=2}^{N-1} a( e_{s(i-1)} , u_{fi}) +  \sum_{i=1}^{N-2} a( e_{s(i+1)} , u_{fi})
\end{equation*}
By definition of the fast interaction term in \eqref{eq:f2e} and inequality $2ab \leq t^{-1} a^2 + t b^2$ for parameter $t>0$ it holds that
\begin{align*}
2\left|a( u_{fi}, e_{s(i-1)} )\right| & \leq t^{-1}\epsilon_f^2 \left\| y_i^{1/2} \nabla u_{fi} \right\|^2_{\Omega_{i}} + t \left\| y_i^{1/2} \nabla e_{s(i-1)} \right\|^2_{\Omega_{i}} \\
2|a( u_{fi}, e_{s(i+1)} )| & \leq t^{-1} \epsilon_f^2 \left\| y_{i+1}^{1/2} \nabla u_{fi} \right\|^2_{\Omega_{i+1}} + t \left\| y_{i+1}^{1/2} \nabla e_{s(i+1)} \right\|^2_{\Omega_{i+1}}
\end{align*}
Thus
\begin{equation*}
\| u_{s} - \tilde{u}_{s} \|^2_E 
\leq
t \sum_{i=1}^{N-1} \| e_{si} \|^2_E + \frac{t^{-1}}{2} \epsilon^2_f \sum_{i=1}^{N-1} \| u_{fi} \|^2_E 
\end{equation*}
By $V_\Gamma$-stability in Lemma~\ref{lemma:Gstab} and slow-fast stability in Lemma~\ref{lemma:stab}, 
\begin{equation*}
\left( 1 - \frac{t}{1-\Lambda} \right) \| u_{s} - \tilde{u}_{s} \|^2_E 
\leq
\frac{t^{-1}}{2} \epsilon^2_f \gamma^2 \| u_\Gamma \|_E^2.
\end{equation*}
Choosing $t = \frac{1}{2}(1-\Lambda)$ completes the proof. 
\end{proof}
\begin{remark} In Lemma~\ref{lemma:serror} we assume that $\epsilon_f$ is sufficiently small. This assumption is required to guarantee that the constants appearing in $V_\Gamma$-stability and slow-fast stability results are positive. The term $\epsilon_f$ is estimated in Section~\ref{sec:sin} and its discrete counterpart in  Section~\ref{sec:SVD}.
\end{remark}

\begin{lemma}[Fast error estimate]
\label{lemma:ferror}
Let $V_{s}$, $V_{si}$, and $V_{fi}$ be as defined in~\eqref{eq:sf_deco} and \eqref{eq:slow_union}. Assume that $V_{si}$ and $V_{fi}$ are $a$-orthogonal for each $i= 1,\ldots, N-1$ and any $y\in D$. Let $u_{fi}(y),\tilde{u}_{fi}(y) \in V_{fi}$ be as defined in \eqref{eq:uexact} and \eqref{eq:tildeuf}, respectively.

\begin{equation*}
\| u_{fi}(y) - \tilde{u}_{fi}(y) \|^2_E \leq  \epsilon^2_f \left( \left\| y_i^{1/2} \nabla u_{i-1}(y) \right\|^2_{\Omega_{i}} + 
\left\| y_{i+1}^{1/2} \nabla u_{i+1}(y) \right\|^2_{\Omega_{i+1}} \right).
\end{equation*}
for any $y\in D$ and $i\in\{1,\ldots,N-1\}$. Here $u_{i}(y) \in V_{\Gamma_i}$ are such that $u_\Gamma(y) = \sum_{i=1}^{N-1} u_i(y)$, $u_0 = u_N = 0$.
\end{lemma} 
\begin{proof} By \eqref{eq:uexact}, \eqref{eq:md:Gproblem}, and orthogonality between $V_{si}$ and $V_{fi}$ it holds that 
\begin{equation*}
\| u_{fi} - \tilde{u}_{fi} \|_E^2 = - a( e_{fi}, u_{i-1}) - a( e_{fi}, u_{i+1}) \quad \mbox{for $e_{fi} = u_{fi} - \tilde{u}_{fi}$}.
\end{equation*}
Here $u_{0} = u_N = 0$. By definition of the fast interaction term in \eqref{eq:f2e} and inequality $2ab \leq a^2 + b^2$,  
\begin{align*}
2\left|a( u_{i-1}, e_{fi} )\right| & \leq \epsilon^2_f \left\| y_i^{1/2} \nabla u_{i-1} \right\|^2_{\Omega_{i}} + \left\| y_i^{1/2} \nabla e_{fi} \right\|^2_{\Omega_{i}} \\
2|a( u_{i+1}, e_{fi} )| & \leq \epsilon^2_f \left\| y_{i+1}^{1/2} \nabla u_{i+1} \right\|^2_{\Omega_{i+1}} + \left\| y_{i+1}^{1/2} \nabla e_{fi} \right\|^2_{\Omega_{i+1}}
\end{align*}
Then 
\begin{equation*}
\| u_{fi} - \tilde{u}_{fi} \|_E^2 \leq 
\frac{1}{2} \| e_{fi} \|_E^2 + \frac{1}{2} \epsilon^2_f \left( \left\| y_i^{1/2} \nabla u_{i-1} \right\|^2_{\Omega_{i}} + 
\left\| y_{i+1}^{1/2} \nabla u_{i+1} \right\|^2_{\Omega_{i+1}} \right)
\end{equation*}
Moving the term $\frac{1}{2} \| e_{fi}\|_E^2$ to LHS completes the proof.
\end{proof}
\begin{theorem}\label{thm:main}  Let $V_{s}$, $V_{si}$, and $V_{fi}$ be as defined in~\eqref{eq:sf_deco} and \eqref{eq:slow_union}. Assume that $V_{si}$ and $V_{fi}$ are $a$-orthogonal for $i\in \{1,\ldots, N-1\}$ and any $y\in D$. In addition let $u_\Gamma \in V_\Gamma$ be the solution to the interface problem and $\tilde{u}_\Gamma$ the approximate interface solution related to the decomposition $V_{s} \oplus V_{f1} \oplus \cdots \oplus V_{f(N-1)}$.  Assume that $\epsilon_f$ is sufficiently small. Then there exists $C$ independent of  $y$, $u$, and $\epsilon_f$ such that
\begin{equation*}
    \| u_\Gamma(y) - \tilde{u}_\Gamma(y) \|_E \leq C \epsilon_f \| u_\Gamma(y) \|_E, 
\end{equation*}
for any $y\in D$. 
\end{theorem}
\begin{proof} By triangle inequality
\begin{equation*}
    \| u_\Gamma - \tilde{u}_\Gamma \|_E \leq \| u_{\Gamma_s} - \tilde{u}_{\Gamma_s} \|_E + \left\| \sum_{i=1}^{N-1} e_{fi} \right\|_E. 
\end{equation*}
By Lemma~\ref{lemma:Gstab} and \ref{lemma:ferror} it holds that 
\begin{equation*}
\left\| \sum_{i=1}^{N-1} e_{fi} \right\|_E^2   
\leq 
(1+\Lambda) \epsilon^2_f \sum_{i=1}^{N-1} \left( \left\| y_i^{1/2} \nabla u_{i-1} \right\|^2_{\Omega_{i}} + 
\left\| y_{i+1}^{1/2} \nabla u_{i+1} \right\|^2_{\Omega_{i+1}} \right)
\end{equation*}
By $V_\Gamma$-stability in Lemma~\ref{lemma:Gstab}, 
\begin{equation*}
\left\| \sum_{i=1}^{N-1} e_{fi} \right\|_E^2   
\leq \frac{1+\Lambda}{1-\Lambda} \epsilon^2_f \| u_\Gamma \|_E^2. 
\end{equation*}
Using Lemma~\ref{lemma:serror} completes the proof. 
\end{proof}

\subsection{Harmonic extension operator}

Theorem~\ref{thm:main} bounds the error between the exact and the approximate interface solutions by the \emph{fast interaction} term
\begin{equation*}
    \epsilon_f = \max_{i=1,\ldots,N-2} \max_{\substack{z_{fi} \in V_{fi} \\ z_{i+1} \in V_{\Gamma_{i+1}} }} \frac{ (\nabla z_{fi},\nabla z_{i+1})_{\Omega_{i+1}} }{ \|\nabla z_{fi} \|_{\Omega_{i+1}} \|\nabla z_{i+1} \|_{\Omega_{i+1}}}.
\end{equation*}
The spaces $V_{fi}$ and $V_{\Gamma_{i+1}}$ are subspaces of $V_\Gamma$ and thus difficult to work with. We alleviate this difficulty by giving a representation of the functions in the space $V_\Gamma$ as harmonic extensions of their interface traces: Let $v_\Gamma \in V_\Gamma$ and define two linear extension operators $Z_{0i} : H^{1/2}_{00}(\Gamma_i) \mapsto H^1(\Omega_{i+1})$ and $Z_{1i} : H^{1/2}_{00}(\Gamma_{i+1}) \mapsto H^1(\Omega_{i+1})$ such that $Z_{0i} h_{i} = w_0$ and $Z_{1i} h_{i+1} = w_1$, where $w_0$ and $w_1$ are weak solutions to problems 
\begin{equation}
\label{eq:Zdef}
\left\{ \begin{aligned} -\Delta w_0 &= 0 & \mbox{in $\Omega_{i+1}$} \\ w_0 &= h_i &\mbox{on $\Gamma_i$} \\ w_0 &= 0 & \mbox{on $\partial \Omega_{i+1} \setminus \Gamma_i$} \end{aligned} \right.
\quad \mbox{and} \quad 
\left\{ \begin{aligned} -\Delta w_1 &= 0 & \mbox{in $\Omega_{i+1}$} \\ w_1 &= h_{i+1} &\mbox{on $\Gamma_{i+1}$} \\ w_1 &= 0 & \mbox{on $\partial \Omega_{i+1} \setminus \Gamma_{i+1}$} \end{aligned} \right.
\end{equation}
Then
\begin{equation*}
v_i|_{\Omega_{i+1}} = Z_{0i} v_i|_{\Gamma_i} \quad \mbox{and} \quad v_{i+1}|_{\Omega_{i+1}} = Z_{1i} v_{i+1}|_{\Gamma_{i+1}}.
\end{equation*}
Recall that $V_{fi}$ is defined as the space of those functions whose trace on $\Gamma_i$ belongs to $Y_{fi}$. Thus, 
\begin{equation*}
\epsilon_f = \max_{i=1,\ldots,N-2} \max_{ \substack{h_{fi} \in Y_{fi} \\ h_{i+1} \in H^{1/2}_{00}(\Gamma_{i+1})} } 
\frac{  (\nabla Z_{0i} h_{fi}, \nabla Z_{1i} h_{i+1} )_{\Omega_{i+1}}   }{ \| \nabla Z_{0i} h_{fi} \|_{\Omega_{i+1}}  \| \nabla Z_{1i} h_{i+1} \|_{\Omega_{i+1}}}.
\end{equation*}

\subsection{Error analysis for sinusoidal slow-fast spaces}
\label{sec:sin}
In this section, we study the error in approximation \eqref{eq:IF_error} for slow-fast decomposition based on an explicit choice of sinusoidal interface spaces. Throughout the section, we assume that each layer $\Omega_{i}$ is a square, i.e, $\Omega_{i} = (i-1,i)\times(0,1)$. The same analysis extends easily to rectangular domains. Let $n_s \in \mathbb{N}$ be a parameter and 
\begin{equation}
\label{eq:cont_if_spaces}
\begin{aligned}
    Y_{fi} &:= \left\{ \sum_{n=n_s+1}^\infty \alpha_n \sin(\pi n x_2) \; | \; \alpha_n \in \mathbb{R} \; \right\} \\  
    Y_{si} &:= \left\{ \sum_{n=1}^{n_s} \alpha_n \sin(\pi n x_2) \; | \; \alpha_n \in \mathbb{R} \; \right\}  
\end{aligned}
\end{equation}
for $i = 1,\ldots,N-1$. Recall that $Y_{fi}$ and $Y_{si}$ define fast interface spaces $V_{fi}$ and the slow space $V_{s}$ via \eqref{eq:sf_deco}. Thus, they also determine the approximate interface solution~\eqref{eq:IF_error}. Error estimate for the approximate interface solution defined by equation~\eqref{eq:cont_if_spaces} follows by bounding the corresponding edge-to-edge interaction and fast interaction terms, $\Lambda$ and $\epsilon_f$, and applying Theorem~\ref{thm:main}. 

First, we give an explicit formula for the harmonic extension operators $Z_{0i}$ and $Z_{1i}$ defined in \eqref{eq:Zdef}. As $Z_{0i}$ and $Z_{1i}$ are linear, it is enough to specify how they operate to $\sin(\pi n x_2)$. By direct computation, it holds that
\begin{align*}
Z_{0i} \sin(\pi n x_2) &= \sin(\pi n x_2) \frac{\sinh(\pi n (i+1-x_1))}{\sinh(\pi n)} \\
Z_{1i} \sin(\pi n x_2) &= \sin(\pi n x_2) \frac{\sinh(\pi n (x_1-i))}{\sinh(\pi n)}.
\end{align*}
Hence, 
\begin{align*}
\nabla Z_{0i} \sin( \pi n x_2) &= 
\frac{\pi n}{\sinh(\pi n)} 
\begin{bmatrix} 
-\sin( \pi n x_2) \cosh(\pi n (i+1-x_1))\\
\cos( \pi n x_2) \sinh(\pi n (i+1-x_1))  
\end{bmatrix}\\
\nabla Z_{1i} \sin( \pi n x_2) &= 
\frac{\pi n}{\sinh(\pi n)} 
\begin{bmatrix} 
\sin( \pi n x_2) \cosh(\pi n (x_1 - i))\\
\cos( \pi n x_2) \sinh(\pi n (x_1 - i)) 
\end{bmatrix}.
\end{align*}
By the above expressions and orthogonality of trigonometric functions it holds that 
\begin{equation*}
(\nabla Z_{0i} h_{si}, \nabla Z_{0i} h_{fi} )_{\Omega_{i+1}} = 0 \quad \mbox{for any $h_{si} \in Y_{si}$ and $h_{fi} \in Y_{fi}$}.
\end{equation*}
and
\begin{equation*}
(\nabla Z_{1i} h_{s(i+1)}, \nabla Z_{1i} h_{f(i+1)} )_{\Omega_{i+1}} = 0 \quad \mbox{for any $h_{s(i+1)} \in Y_{s(i+1)}$ and $h_{f(i+1)} \in Y_{f(i+1)}$}.
\end{equation*}
Hence, the spaces $V_{si}$ and $V_{fi}$ corresponding to $Y_{si}$ and $Y_{fi}$ in \eqref{eq:cont_if_spaces} are $a$-orthogonal for any $y\in D$. We proceed to estimate the fast interaction term

\begin{equation*}
\epsilon_f = \max_{ \substack{h_{fi} \in Y_{fi} \\ h_{i+1} \in H^{1/2}_{00}(\Gamma_{i+1})} } 
 \frac{  (\nabla Z_{0i} h_{fi}, \nabla Z_{1i} h_{i+1} )_{\Omega_{i+1}}   }{ \| \nabla Z_{0i} h_{fi} \|_{\Omega_{i+1}}  \| \nabla Z_{1i} h_{i+1} \|_{\Omega_{i+1}}}.
\end{equation*}
Note that the fast interaction term is the same for any $i \in\{1,\ldots,N-2\}$. By direct computation it holds that 
\begin{align*}
\| \nabla Z_{0i} h_{fi} \|_{\Omega_{i+1}}^2 &= \sum_{n=n_s+1}^\infty \alpha_n^2 \| \nabla Z_{0i} \sin(\pi n x_2) \|_{\Omega_{i+1}}^2= \sum_{n=n_s+1}^\infty \alpha_n^2 \frac{\pi n \sinh(2\pi n)}{4\sinh^2(\pi n)} \\
\| \nabla Z_{1i} h_{i+1} \|_{\Omega_{i+1}}^2 &= \sum_{n=1}^\infty \beta_n^2 \| \nabla Z_{1i} \sin(\pi n x_2) \|_{\Omega_{i+1}}^2= \sum_{n=1}^\infty \beta_n^2 \frac{\pi n \sinh(2\pi n)}{4\sinh^2(\pi n)} \\
(\nabla Z_{0i} h_{fi}, \nabla Z_{1i} h_{i+1} )_{\Omega_{i+1}} &= \sum_{n=n_s+1}^\infty \alpha_n \beta_n (\nabla Z_{0i} \sin(\pi n x_2),\nabla Z_{1i} \sin(\pi n x_2))_{\Omega_{i+1}}\\
&= -\sum_{n=n_s+1}^\infty \alpha_n \beta_n \frac{\pi n}{2\sinh(\pi n)}
\end{align*}
for any $h_{fi} \in Y_{fi}$ and $h_{i+1} \in H^{1/2}_{00}(\Gamma_{i+1})$ such that 
\begin{equation*}
h_{fi} = \sum_{n=n_s+1}^\infty \alpha_n \sin(\pi n x_2)
\quad \mbox{and} \quad 
h_{i+1} = \sum_{n=1}^\infty \beta_n \sin(\pi n x_2).
\end{equation*}
%
%
%
By utilizing the above formulas, we obtain the following estimate for the interaction term. 

\begin{lemma} \label{lemma:exponential} Assume that $\Omega_{i} = (i-1,i)\times(0,1)$ for $i=1,\ldots,N$. Let $n_s \in \mathbb{N}$, spaces $Y_{si},Y_{fi}$ be as defined in \eqref{eq:cont_if_spaces}, $\epsilon_f \equiv \epsilon_f(n_s)$ the fast interaction term defined in \eqref{eq:f2e}, and $\Lambda$ the edge-to-edge interaction term defined in \eqref{eq:e2e}. Then there holds that
\begin{equation*}
\epsilon_f(n_s) \leq C e^{-\pi (n_s+1)} \quad  \mbox{and} \quad \Lambda \leq \frac{2 e^{-\pi}}{1-e^{-2\pi}}.
\end{equation*}
\end{lemma}
\begin{proof} 
First, we note the lower limit 
\begin{equation*}
    \frac{\sinh(2\pi n)}{\sinh^2(\pi n)} \geq 2.
\end{equation*}
This allows us to estimate using the Cauchy-Schwarz inequality:
\begin{align*}
    \frac{|(\nabla Z_{0i} h_{fi},\nabla Z_{1i} h_{1+1} )_{\Omega_{i+1}}|}{\| \nabla Z_{0i} h_{fi} \| \| \nabla Z_{1i} h_{i+1} \|}
    &\leq \frac{ \sum_{n=n_s+1}^{\infty} \frac{\alpha_n \beta_n n}{\sinh(\pi n)}}{\left( \sum_{n=n_s+1}^{\infty} \alpha_n^2 n \right)^{1/2} \left( \sum_{n=n_s+1}^{\infty} \beta_n^2 n \right) ^{1/2}}\\
    &\leq \frac{\frac{1}{\sinh(\pi (n_s+1))} \left( \sum_{n=n_s+1}^{\infty} \alpha_n^2 n \right)^{1/2} \left( \sum_{n=n_s+1}^{\infty} \beta_n^2 n  \right)^{1/2}}{\left( \sum_{n=n_s+1}^{\infty} \alpha_n^2 n \right)^{1/2} \left( \sum_{n=n_s+1}^{\infty} \beta_n^2 n \right) ^{1/2}}\\
    &= \frac{1}{\sinh(\pi (n_s+1))}
    \leq C e^{-\pi (n_s+1)}
\end{align*}
where the constant can be taken as $C=2/(1-e^{-2\pi}) \approx 2.004$. Choosing $n_s=0$ corresponds to the case $Y_{fi} = H^{1/2}_{00}(\Gamma_i)$ and gives the estimate for $\Lambda$. 
\end{proof}
The above results yield the following theorem:
\begin{theorem} \label{thm:sine_error} Assume that $\Omega_{i} = (i-1,i)\times(0,1)$ for $i=1,\ldots,N$. Let $n_s\in \mathbb{N}$, the spaces $Y_{si}$ and $Y_{fi}$ be as defined in \eqref{eq:cont_if_spaces}, $V_{si}$ and  $V_{fi}$ as in \eqref{eq:sf_deco}, and $V_s$ as in \eqref{eq:slow_union}. In addition, let $\tilde{u}_\Gamma(y) \in V_\Gamma$ be the approximate interface solution \eqref{eq:IF_error}. Then there exists $C>0$ independent of $n_s$, $y$, $f$ and $u$ s.t.
\begin{equation*}
 \| u_\Gamma(y) - \tilde{u}_\Gamma(y) \|_E \leq C e^{-\pi (n_s+1)} \| u_\Gamma(y) \|_E
\end{equation*}
for any $y \in D$.
\end{theorem}
\begin{proof} First, we verify assumptions made in Theorem~\ref{thm:main}. By direct computation, the spaces $V_{si}$ and $V_{fi}$ are $a$-orthogonal. The assumption on $\epsilon_f$ being sufficiently small in Theorem~\ref{thm:main} was made to guarantee that the stability constant $\gamma$ appearing in Lemma~\ref{lemma:stab} is positive. In the current case this is always true.    

The result follows by combining Theorem~\ref{thm:main} and Lemma~\ref{lemma:exponential}. 
\end{proof}

\section{Finite element realisation}
\label{sec:fem}

In this section, we construct approximate parameter-to-solution maps in the context of the finite element method, see e.g., \cite{Br:2001,BrSc:2008}. The finite element approximation to \eqref{eq:prob} is: find $u_h:D\mapsto V_h$ such that $u_h(y) \in V_h$ satisfies  
\begin{equation}
\label{eq:fem_prob}
    a(\nabla u_h(y), \nabla v_h) = (f,v_h) \quad \mbox{for all $v_h \in V_h$ and $y\in D$}.  
\end{equation}
Here $V_h \subset H^1_0(\Omega)$, $V_h = \mathop{span}_{i=1}^n \{\varphi_i\}$, is the finite element space. We assume that $V_h$ conforms to interfaces $\Gamma_i$ and is locally mirror symmetric in the sense of Assumption~\ref{ass:symm}. 
\begin{assumption} 
\label{ass:symm}
Let $v_h \in V_h$. Then there exists $w_h \in V_h$ such that
\begin{equation*}
    w_h = v_h \quad \mbox{in $\Omega_i$} \quad \mbox{and} \quad w_h(i-t,x_2) = w_h(i+t,x_2)
\end{equation*}
for any $t,x_2$ such that $(i-t, x_2) \in \Omega_i$, $i=1,...,N-1$. 
\end{assumption}
\begin{figure}
    \centering
    \includegraphics[width=0.5\textwidth]{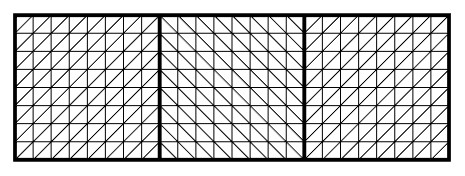}
    \caption{A locally mirror symmetric finite element mesh for the $3\times 1$ square grid. Standard uniform order finite element spaces associated to this mesh satisfy Assumption~\ref{ass:symm}.}
    \label{fig:symm_mesh}
\end{figure}
\noindent Assumption~\ref{ass:symm} is satisfied if $V_h$ is a finite element space of uniform degree associated to a mesh that is locally mirror symmetric with respect to every interface, see Figure \ref{fig:symm_mesh}. Note that we are still assuming for simplicity that every interface is one unit away from the previous one. However, this is not necessary and could be relaxed. 

Let $v_h \in V_h$ and $\vec{\alpha} \in \mathbb{R}^n$.  We write 
\begin{equation*}
\vec{\alpha} \sim v_h \quad \mbox{if} \quad v_h = \sum_{i=1}^n \alpha_i \varphi_i.
\end{equation*} 
This is, $\vec{\alpha} \in \mathbb{R}^n$ is the coordinate vector of $v_h \in V_h$ in the basis $\{ \varphi_1,\ldots,\varphi_n\}$. By direct computation, for any bilinear form $b:V_h \times V_h \mapsto \mathbb{R}$ there exists a matrix $\mathsf{B} \in \mathbb{R}^{n\times n}$ such that
\begin{equation*}
    \vec{\alpha}^T \mathsf{B} \vec{\beta} = b(w_h, z_h) \quad \mbox{for any $w_h, z_h \in V_h$ and $\vec{\alpha} \sim w_h$, $\vec{\beta} \sim z_h$}. 
\end{equation*}

We proceed by repeating key definitions and results of Section~\ref{sec:error_estimates} in the discrete setting. Due to Assumption~\ref{ass:symm}, Lemma \ref{lemma:average} and Theorem \ref{thm:2x1} hold. The discrete subdomain-interface decomposition is $V_h = V_{h,\Omega_1}\oplus \cdots  \oplus V_{h,\Omega_N} \oplus V_{h,\Gamma}$, where 
\begin{equation*}
V_{h,\Omega_{i}} =  \{\; v_h \in V_h \; |\; supp(v_h) \subset \Omega_{i} \; \} 
\end{equation*}
and
\begin{equation*}
V_{h,\Gamma} = \{\; v_h \in V_h \; |\; (\nabla v_h,\nabla v_{h,\Omega_i})_\Omega = 0 \quad \mbox{for any $v_{h,\Omega_i} \in V_{h,\Omega_i}$ and $i=1,\ldots,N$} \; \}. 
\end{equation*}
Observe that $V_{h,\Gamma} \not\subset V_\Gamma$. Clearly, the solution to \eqref{eq:fem_prob} satisfies 
\begin{equation}
    \label{eq:disc_sub_int_deco}
    u_h(y) = \sum_{i=1}^N y_i^{-1} w_{h,\Omega_i} + u_{h,\Gamma}(y) \quad \mbox{for any $y\in D$}.
\end{equation}
Here $w_{h,\Omega_{i}} \in V_{h,\Omega_i}$ satisfies $(\nabla w_{h,\Omega_i},\nabla v_{h,\Omega_i} )_\Omega = (f,v_{h,\Omega_i} )$ for each $v_{h,\Omega_i} \in V_{h,\Omega_i}$ and the discrete interface solution $u_{h,\Gamma}(y)\in V_{h,\Gamma}$ satisfies \eqref{eq:fem_prob} posed in $V_{h,\Gamma}$. 

Identical to the continuous case, we decompose the discrete interface space as $V_{h,\Gamma} = V_{h,s} \oplus V_{h,f1} \oplus \cdots \oplus V_{h,f(N-1)}$. Let the discrete slow and fast interface spaces $Y_{h,si}, Y_{h,fi}$ satisfy $V_h|_{\Gamma_i} = Y_{h,si} \oplus Y_{h,fi}$ and define
\begin{equation}
\label{eq:disc_sf}
\begin{aligned}
V_{h,fi} &= \{ v \in V_{h,\Gamma} \; | \; v|_{\Gamma_i} \in Y_{h,fi} \mbox{ and } v|_{\Gamma_j} = 0 \mbox{ for $j\neq i$} \; \} \\
V_{h,si} &= \{ v \in V_{h,\Gamma} \; | \; v|_{\Gamma_i} \in Y_{h,si} \mbox{ and } v|_{\Gamma_j} = 0 \mbox{ for $j\neq i$} \; \} \\
V_{h, s} &= \bigoplus_{i=1}^{N-1} V_{h,si}.
\end{aligned}    
\end{equation}
The approximate interface solution related to the above discrete slow-fast decomposition is 
\begin{equation}
\label{eq:disc_approx}
\tilde{u}_{h,\Gamma}(y) = \sum_{i=1}^{N-1} \frac{2}{y_i+y_{i+1}} w_{h,fi} + \tilde{u}_{h,s}(y), 
\end{equation}
where $\tilde{u}_{h,s}(y) \in V_{h,s}$ satisfies 
\begin{equation*}
a(\tilde{u}_{h,s}(y),v_{h,s} ) = (f,v_{h,s}) \quad \mbox{for any $v_{h,s} \in V_{h,s}$ and $y\in D$}.
\end{equation*}
Each $w_{h,fi} \in V_{h,fi}$ in \eqref{eq:disc_approx} satisfies 
\begin{equation*}
(\nabla w_{h,fi},\nabla v_{h,fi} )_\Omega = (f,v_{h,fi})_\Omega \quad \mbox{for any $v_{h,fi} \in V_{h,fi}$ and $y\in D$}. 
\end{equation*}

The error analysis given in Section~\ref{sec:error_estimates} applies to the finite dimensional case with the exception that the interaction terms $\Lambda$ and $\epsilon_f$ are replaced by their discrete counterparts: 
\begin{equation*}
    \Lambda_h = \max_{i=1,\ldots,N-2} \max_{\substack{z_{i} \in V_{h,\Gamma_i} \\ z_{i+1} \in V_{h,\Gamma_{i+1}}}} \frac{ (\nabla z_i,\nabla z_{i+1})_{\Omega_{i+1}} }{ \|\nabla z_{i} \|_{\Omega_{i+1}} \|\nabla z_{i+1} \|_{\Omega_{i+1}}}
\end{equation*}
and
\begin{equation}
\label{eq:discrete_f2e}
    \epsilon_{h,f} = \max_{i=1,\ldots,N-2} \max_{\substack{z_{fi} \in V_{h,fi} \\ z_{i+1} \in V_{h,\Gamma_{i+1} }}} \frac{ (\nabla z_{fi},\nabla z_{i+1})_{\Omega_{i+1}} }{ \|\nabla z_{fi} \|_{\Omega_{i+1}} \|\nabla z_{i+1} \|_{\Omega_{i+1}}}.
\end{equation}
Here $V_{h,\Gamma_{i}} := \{ \;  v \in V_{h,\Gamma} \; | \; v|_{\Gamma_j} = 0 \mbox{ for $j\neq i$} \; \}$. Error estimate given in Theorem~\ref{thm:main} is valid in the finite element setting, if $V_{h,si}$ and $V_{h,fi}$ are $a$-orthogonal, $\Lambda_h < 1$ and $\epsilon_{h,f}$ is sufficiently small. Then there exists $C$ such that
\begin{equation*}
\| u_{h,\Gamma}(y) - \tilde{u}_{h,\Gamma}(y) \|_E \leq C \epsilon_{h,f} \| u_{h,\Gamma}(y) \|_E \quad \mbox{for any $y\in D$},
\end{equation*}
where $\tilde{u}_{h,\Gamma}(y) \in V_{h,\Gamma}$ is the approximate discrete interface solution defined in~\eqref{eq:disc_approx}.

\subsection{Discrete slow-fast spaces}
\label{sec:SVD}
In this section, we give a method for computing discrete slow-and-fast spaces that minimize the fast interaction term in \eqref{eq:discrete_f2e} for a given dimension of the slow space. We represent functions in $V_{h,i}$ and $V_{h,i+1}$ using discrete harmonic extension operators $Z_{h,0i} : V_h|_{\Gamma_{i}} \mapsto V_h|_{\Omega_{i+1}}$ and $Z_{h,1i} : V_h|_{\Gamma_{i+1}} \mapsto V_h|_{\Omega_{i+1}}$. The operator $Z_{h,0i}$ is defined as $Z_{h,0i} g_h = w_h$, where 
\begin{equation*}
w_h \in \{ \; v_h \in V_h|_{\Omega_{i+1}} \; |  \; \mbox{$v_h|_{\Gamma_i} = g_h$ and $v_h|_{\Gamma_{i+1}} = 0$} \; \}
\end{equation*}
satisfies 
\begin{equation*}
(\nabla w_h, \nabla v_{h,\Omega_{i+1}})_{\Omega_{i+1}} = 0 \quad \mbox{for any $v_{h,\Omega_{i+1}} \in V_{h,\Omega_{i+1}}$}. 
\end{equation*}
This is, $Z_{h,0i} g_h$ is the discrete harmonic extension of $g_h$ to $\Omega_{i+1}$. The definition of $Z_{h,1i}$ is similar and not repeated. 

The fast interaction term is formulated using matrix notation. Let matrices $\mathsf{C}^{(i)}_{nm}$ for $i = 1,\ldots,N-1$ and $n,m \in \{0,1\}$ satisfy
\begin{equation}
\label{eq:def:Cnm}
 \vec{\alpha}^T \mathsf{C}^{(i)}_{nm} \vec{\beta} = (\nabla Z_{h,ni} g_{i+n} ,\nabla Z_{h,mi} g_{i+m} )_{\Omega_{i+1}}
\end{equation}
for any $g_{i+n} \in V_h|_{\Gamma_{i+n}}$, $g_{i+m} \in V_h|_{\Gamma_{i+m}}$ and $\vec{\alpha} \sim g_{i+n}$, $\vec{\beta} \sim g_{i+m}$. In the rest of this section, we omit the superscript $(i)$ when possible.  This is, we set $C_{nm} \equiv C_{nm}^{(i)}$ for $n,m\in\{0,1\}$.  

Let $\mathsf{R}_0$ and $\mathsf{R}_1$ be the Cholesky factors of $\mathsf{C}_{00}$ and $\mathsf{C}_{11}$, respectively. Then 
\begin{equation*}
\frac{(\nabla Z_{h,0i} g_{i},\nabla Z_{h,1i} g_{i+1} )_{\Omega_{i+1}}}{\| \nabla Z_{h,0i} g_{i} \| \| \nabla Z_{h,1i} g_{i+1} \|} = 
\frac{\vec{\alpha}^T \mathsf{C}_{01} \vec{\beta} }{ | \mathsf{R}_{0} \vec{\alpha} | | \mathsf{R}_{1} \vec{\beta} |}. 
\end{equation*}
For any $g_i \in V_{h}|_{\Gamma_{i}}, g_{i+1} \in V_{h}|_{\Gamma_{i+1}}$ and $\vec{\alpha} \sim g_{i}$, $\vec{\beta} \sim g_{i+1}$. Let $\mathsf{U} \mathsf{\Sigma} \mathsf{V}^T$ be the SVD of $\mathsf{R}_0^{-T} \mathsf{C}_{01} \mathsf{R}_1^{-1}$. Then: 
\begin{equation*}
\frac{(\nabla Z_{h,0i} g_{i},\nabla Z_{h,1i} g_{i+1} )_{\Omega_{i+1}}}{\| \nabla Z_{h,0i} g_{i} \| \| \nabla Z_{h,1i} g_{i+1} \|} =
\frac{\vec{\alpha}^T  \mathsf{R}_{0}^T \mathsf{U} \mathsf{\Sigma} \mathsf{V}^T \mathsf{R}_{1} \vec{\beta} }{ | \mathsf{R}_{0} \vec{\alpha} | | \mathsf{R}_{1} \vec{\beta} |}.
\end{equation*}
We use this expression to choose the slow and fast interface spaces $Y_{h,fi},Y_{h,si}$ that define the discrete slow-fast decomposition via \eqref{eq:disc_sf}. Let 
$i = 1,\ldots,N-1$ and $r_i \in \mathbb{N}$ be the truncation index. Define  
\begin{equation}
\label{eq:disc_if_spaces}
\begin{aligned}
Y_{h,si} & := \{ \; g_{h,si} \in V_h|_{\Gamma_i} \; | \; g_{h,si} \sim \vec{\alpha} \mbox{ for } \mathsf{R}_0^{(i)} \vec{\alpha} \in \mathop{range}( \mathsf{U}^{(i)}(:,1:r_i) \; \} \\
Y_{h,fi} & := \{ \; g_{h,fi} \in V_h|_{\Gamma_i} \; | \; g_{h,fi} \sim \vec{\alpha} \mbox{ for } \mathsf{R}_0^{(i)} \vec{\alpha} \in \mathop{range}( \mathsf{U}^{(i)}(:,(r_i+1):end)) \; \}.
\end{aligned}
\end{equation}
Recall that matrices with superscript $(i)$ correspond to $C_{nm}^{(i)}$. 

We proceed to show that $V_{h,fi}$ and $V_{h,si}$ related to these interface spaces are $a$-orthogonal, interaction terms satisfy $\Lambda_h = \sigma_1$, $\epsilon_{h,f} = \sigma_{r_i+1}$, and these interface spaces have the smallest fast interaction term with the given slow space dimension.
\begin{lemma} Let the interface spaces $Y_{h,si}$ and $Y_{h,fi}$ be as defined in \eqref{eq:disc_if_spaces}. Then the corresponding spaces $V_{h,si}$ and $V_{h,fi}$ defined in \eqref{eq:disc_sf} are $a$ - orthogonal for any $y\in D$. 
\end{lemma}
\begin{proof} Let $g_{si} \in V_{h,si}$ and $g_{fi} \in V_{h,fi}$. By symmetry  
\begin{equation*}
    a(g_{si},g_{fi}) = (y_{i} + y_{i+1}) (\nabla Z_{h,0i} z_{fi},\nabla Z_{h,0i} z_{si})_{\Omega_{i+1}}
\end{equation*}
for some $z_{fi} \in Y_{h,fi}$ and $z_{si} \in Y_{h,si}$. By \eqref{eq:disc_if_spaces}
\begin{equation*}
z_{si} \sim \mathsf{R}_0^{-1} \mathsf{U}(:,1:r_i) \vec{\alpha} \quad \mbox{and} \quad z_{fi} \sim \mathsf{R}_0^{-1} \mathsf{U}(:,(r_i+1):end) \vec{\beta}
\end{equation*}
for some $\vec{\alpha},\vec{\beta}$ of appropriate dimension. Using \eqref{eq:def:Cnm},  gives 
\begin{equation*}
(\nabla Z_{0i} z_{si}, \nabla Z_{0i} z_{fi} ) = \vec{\alpha}^T \mathsf{U}(:,1:r_i)^T \mathsf{R}_0^{-T} \mathsf{C}_{00} \mathsf{R}_0^{-1} \mathsf{U}(:,(r_i+1):end) \vec{\beta}
\end{equation*}
Recalling that $\mathsf{C}_{00}=\mathsf{R}_0^T \mathsf{R}_0$ and $\mathsf{U}^T \mathsf{U} = \mathsf{I}$ completes the proof. 
\end{proof}

\begin{lemma} Let the interface spaces $Y_{h,si}$ and $Y_{h,fi}$ be as defined in \eqref{eq:disc_if_spaces}. Then the discrete fast interaction term related to the $i$th interface satisfies 
\begin{equation*}
\max_{\substack{ z_{fi} \in Y_{h, fi} \\ z_{i+1} \in V_h|_{\Gamma_{i+1}} } } \frac{ (\nabla Z_{h,0i} z_{fi},\nabla Z_{h,1i} z_{i+1} )_{\Omega_{i+1}}}{\| \nabla Z_{h,0i} z_{fi} \|_{\Omega_{i+1}} \| \nabla Z_{h,1i} z_{i+1} \|_{\Omega_{i+1}}} = \sigma^{(i)}_{r_i+1}.    
\end{equation*}
Here $\{ \sigma_k^{(i)} \}$  are the singular values of $(\mathsf{R}_0^{(i)})^{-T} \mathsf{C}^{(i)}_{01} (\mathsf{R}^{(i)}_1)^{-1}$.
\end{lemma}
\begin{proof}  
Let $z_{fi} \in Y_{h,fi}$ and $z_{i+1} \in V_{i}|_{\Gamma_{i+1}}$. Then there holds that
\begin{equation*}
z_{fi} \sim \mathsf{R}_0^{-1} \mathsf{U}(:,(r_i+1):end) \vec{\alpha} \quad \mbox{and} \quad z_{i+1} \sim \mathsf{R}_1^{-1} \vec{\beta} 
\end{equation*}
for some $\vec{\alpha},\vec{\beta}$. By \eqref{eq:def:Cnm}, 
\begin{equation*}
\begin{aligned}
(\nabla Z_{h,0i} z_{fi},\nabla Z_{h,1i} z_{i+1} ) & = \vec{\alpha}^T \mathsf{U}(:,(r_i+1):end)^T \mathsf{R}_0^{-T} \mathsf{C}_{01} \mathsf{R}_1^{-1} \vec{\beta} \\
\| \nabla Z_{h,0i} z_{fi}\|_{\Omega_{i+1} } & = |\vec{\alpha}| \\ 
\| \nabla Z_{h,1i} z_{i+1}\|_{\Omega_{i+1} } & = |\vec{\beta}| \\ 
\end{aligned}
\end{equation*}
Recall that $\mathsf{U} \mathsf{\Sigma} \mathsf{V}^T = \mathsf{R}_0^{-T} \mathsf{C}_{01} \mathsf{R}_1^{-1}$ so that 
\begin{equation*}
(\nabla Z_{h,0i} z_{fi},\nabla Z_{h,1i} z_{i+1} ) = \vec{\alpha}^T \mathsf{\Sigma}((r_i+1):end,(r_i+1):end) \mathsf{V}(:,(r_i+1):end)^T \vec{\beta}.
\end{equation*}
By Cauchy-Schwarz inequality, $|(\nabla Z_{h,0i} z_{fi},\nabla Z_{h,1i} z_{i+1} )| \leq \sigma_{r_i+1} | \vec{\alpha} | | \vec{\beta} |$. This estimate becomes equality by choosing $\vec{\alpha} = \vec{e}_1$ and $\vec{\beta} = \mathsf{V}(:,r_i+1)$. This completes the proof. 
\end{proof}
Observe that $\Lambda_h = \sigma_1$. Hence, the largest singular value corresponds to the edge-to-edge interaction term. It is straightforward to see the spaces in~\eqref{eq:disc_if_spaces} are the best in the following sense. 
\begin{lemma} Let the interface spaces $Y_{h,si}$ and $Y_{h,fi}$ be as defined in \eqref{eq:disc_if_spaces}. Then the discrete fast interaction term related to the $i$th interface satisfies
\begin{equation*}
\max_{\substack{ z_{fi} \in X \\ z_{i+1} \in V_h|_{\Gamma_{i+1}} } } \frac{ (\nabla Z_{h,0i} z_{fi},\nabla Z_{h,1i} z_{i+1} )_{\Omega_{i+1}}}{\| \nabla Z_{h,0i} z_{fi} \|_{\Omega_{i+1}} \| \nabla Z_{h,1i} z_{i+1} \|_{\Omega_{i+1}}} \geq \sigma_{r_i+1}
\end{equation*}
for any $X \subset V_{h}|_{\Gamma_i}$ s.t. $\mathop{dim}(X) = \mathop{dim}(Y_{h,fi})$. Here $\{ \sigma_k^{(i)} \}$  are the singular values of $(\mathsf{R}_0^{(i)})^{-T} \mathsf{C}^{(i)}_{01} (\mathsf{R}^{(i)}_1)^{-1}$.  
\end{lemma}
\begin{proof}
Let $z_{fi} \in X$ and $z_{i+1} \in V_{i}|_{\Gamma_{i+1}}$. Then there holds that
\begin{equation*}
z_{fi} \sim \mathsf{R}_0^{-1} \mathsf{Q} \vec{\alpha} \quad \mbox{and} \quad z_{i} \sim \mathsf{R}_1^{-1} \vec{\beta} 
\end{equation*}
for some unitary matrix $\mathsf{Q}$ and some $\vec{\alpha},\vec{\beta}$ of appropriate dimension. By \eqref{eq:def:Cnm}, 
\begin{equation*}
\begin{aligned}
(\nabla Z_{h,0i} z_{fi},\nabla Z_{h,1i} z_{i+1} ) & = \vec{\alpha}^T \mathsf{Q}^T \mathsf{R}_0^{-T} \mathsf{C}_{01} \mathsf{R}_1^{-1} \vec{\beta} \\
\| \nabla Z_{h,0i} z_{fi}\|_{\Omega_{i+1} } & = |\vec{\alpha}| \\ 
\| \nabla Z_{h,1i} z_{i+1}\|_{\Omega_{i+1} } & = |\vec{\beta}|.
\end{aligned}
\end{equation*}
Then the discrete fast interaction term related to the $i$th interface is 
\begin{equation*}
\max_{\substack{ z_{fi} \in X \\ z_{i+1} \in V_h|_{\Gamma_{i+1}} } } \frac{ (\nabla Z_{h,0i} z_{fi},\nabla Z_{h,1i} z_{i+1} )_{\Omega_{i+1}}}{\| \nabla Z_{h,0i} z_{fi} \|_{\Omega_{i+1}} \| \nabla Z_{h,1i} z_{i+1} \|_{\Omega_{i+1}}} = \max_{\vec{\alpha},\vec{\beta}} \frac{ \vec{\alpha}^T \mathsf{Q}^T \mathsf{R}_0^{-T} \mathsf{C}_{01} \mathsf{R}_1^{-1} \vec{\beta} }{|\vec{\alpha}||\vec{\beta}|}    
\end{equation*}
and further 
\begin{equation*}
\max_{\vec{\alpha},\vec{\beta}} \frac{ \vec{\alpha}^T \mathsf{Q}^T \mathsf{R}_0^{-T} \mathsf{C}_{01} \mathsf{R}_1^{-1} \vec{\beta} }{|\vec{\alpha}||\vec{\beta}|} = \left( \max_{\vec{x} \in \mathop{range}(\mathsf{Q})} \frac{ \vec{x}^T \mathsf{R}_0^{-T} \mathsf{C}_{01} \mathsf{R}_1^{-1} \mathsf{R}_1^{-T} \mathsf{C}_{01}^T \mathsf{R}_0^{-1} \vec{x} }{ \vec{x}^T \vec{x} } \right)^{1/2}. 
\end{equation*}
Recalling that $\sigma^2_{r_i+1}$ is the $(r_i+1)$th largest eigenvalue of $\mathsf{R}_0^{-T} \mathsf{C}_{01} \mathsf{R}_1^{-1} \mathsf{R}_1^{-T} \mathsf{C}_{01}^T \mathsf{R}_0^{-1}$, $rank(Q) = r_i+1$, and applying the Courant-Fisher min-max theorem completes the proof. 
\end{proof}

The above results are collected to the following Theorem.

\begin{theorem} \label{thm:svd_error} Let the spaces $Y_{h,si}$ and $Y_{h,fi}$ be as defined in \eqref{eq:disc_if_spaces}, and $V_{h,si}$, $V_{h,fi}$ and $V_{h,s}$ as in \eqref{eq:disc_sf}. In addition, let $\tilde{u}_{h,\Gamma } : D \mapsto  V_{h,\Gamma}$ be the approximate discrete interface solution defined in \eqref{eq:disc_approx}, $\sigma_1 = max_{i=1,\ldots,N-1} \sigma^{(i)}_{1}$, and $\sigma_{r+1} = max_{i=1,\ldots,N-1} \sigma^{(i)}_{r_i+1}$. Assume that $\sigma_1 < 1$ and $\sigma_{r+1}$ is sufficiently small.  Then there exits $C$ such that
\begin{equation*}
 \| u_{h,\Gamma}(y) - \tilde{u}_{h,\Gamma}(y) \|_E \leq C \sigma_{r+1} \| u_{h,\Gamma}(y) \|_E \quad \mbox{for any $y \in D$.}
\end{equation*}
 Here $\{ \sigma_k^{(i)} \}$  are the singular values of $(\mathsf{R}_0^{(i)})^{-T} \mathsf{C}^{(i)}_{01} (\mathsf{R}^{(i)}_1)^{-1}$.
\end{theorem}
\begin{proof} This theorem follows by noting that $\sigma_1$ and $\sigma_{r+1}$ correspond to $\Lambda_h$ and $\epsilon_{h,f}$, respectively, and using the discrete version of Theorem~\ref{thm:main}.
\end{proof}

\section{Applications}

In this section we apply our results to estimating the Kolmogorov $n$-width. Estimates for the Kolmogorov $n$-width are important as they are used to analyse the efficiency of reduced basis methods for parametric PDEs, see \cite{BaCo:17}. 

Numerical solution methods for parametric PDEs, such as \eqref{eq:prob}, seek for an approximate parameter-to-solution map $u_n : D \mapsto H^1_0(\Omega)$ of the form 
\begin{equation*}
    u_n(y) = \sum_{k=1}^n v_k(x) \phi_k(y).
\end{equation*}
The sets of functions $\{v_k\}$ and $\{\phi_k\}$ depend on the applied method. For example, in \emph{sparse collocation methods} $\{ v_k \}$ are the evaluations of the solution at sparse collocation points and $\{ \phi_k \}$ are the corresponding interpolation functions, see \cite{NoTeWe:08}. In \emph{reduced basis methods}, $\{v_k\}$ are obtained by sampling $u$ at selected points, e.g. by utilizing a greedy method, see \cite{QuMaNe:16}. The weights $\{ \phi_k \}$ are then computed pointwise by solving \eqref{eq:prob} in the subspace $V_n := \mathop{span} \{ v_1,\ldots,v_n \}$.

The efficiency of reduced basis methods depends on the Kolmogorov n-width \cite{Pi:1985}, 
\begin{equation*}
    d_n := \min_{ \mathop{dim}(X) = n} \max_{y\in D} \| u(y) - P_X u(y) \|,  
\end{equation*}
where $\| \cdot \|$ is a suitable norm induced by an inner product and $P_X$ is the orthogonal projection to the subspace $X$ in that inner product. Here we use a relative version in the parameter dependent energy norm and define 
\begin{equation*}
    \hat{d}_n := \min_{ \mathop{dim}(X) = n} \max_{y\in D} \frac{\| u(y) - P_E u(y) \|_E}{\|u(y)\|_E}.   
\end{equation*}
Here $P_E$ is the orthogonal projection to $X$ in the parameter dependent $a(\cdot,\cdot)$-inner product. 

Kolmogorov $n$-width measures the approximation error produced by a reduced basis method utilizing an optimal $n$-dimensional subspace $X$. The greedy sampling process can produce subspaces that yield an error comparable to $d_n$, see \cite{BiCo:11,DePePr:12}. Hence, it is interesting to estimate $d_n$ as it also describes the optimal accuracy achievable by a reduced basis method using an $n$-dimensional method subspace. 

In literature, Kolmogorov $n$-width is typically estimated by constructing a subspace $X$ utilizing the functions $v_k$ that are related to polynomial approximation of $u(y)$ \cite{CoDeSc:2010, CoDeSc:2011}. Naturally, $d_n$ is bounded from above with the associated polynomial approximation error. However, the results in \cite{BaCoDa:2018} indicate that in the case of piecewise constant parameters in square grids subspace-based approximations outperform polynomial approximations. The Kolmogorov $n$-width related to parametric PDEs with affine parameter dependency is studied in \cite{BaCo:17} by a truncated Neumann series approximation and analyzing the truncation error.

It immediately follows from Theorem \ref{thm:2x1} that for two mirror symmetric layers $\hat{d}_3=0$. For $N$ layers, the Theorems \ref{thm:sine_error} and \ref{thm:svd_error} immediately give bounds for the relative Kolmogorov $n$-width of problem \eqref{eq:prob} in the continuous space for rectangular layers and in the finite element space for any mirror symmetric layers, respectively.

\begin{theorem} Consider the problem \eqref{eq:prob} with rectangular layers. There exists a constant $C$ such that 
\begin{equation*}
    \hat{d}_{ 2N - 1 + (N-1) n_s} \leq  C e^{-\pi (n_s+1)}
\end{equation*}
for any $n_s \in \mathbb{N}$.
\end{theorem}
\begin{proof} The result follows from Equations \eqref{eq:md:sub_int_deco}, \eqref{eq:md:sub_sol} and Theorem \ref{thm:sine_error} by choosing the subspace
\begin{equation*}
X = \mathop{span} \{ w_{\Omega_1},\ldots,w_{\Omega_N}, w_{f1},\ldots,w_{f(N-1)}, \tilde{u}_s(y) \}.
\end{equation*}
Then the dimension of this space is $\dim(X) = N + (N-1) + (N-1) n_s$.
\end{proof}

\begin{theorem} Consider the problem \eqref{eq:fem_prob}. Make the same assumptions and use the same notation as in Theorem \ref{thm:svd_error}. Then there exists a constant $C$ such that

\begin{equation*}
    \hat{d}_{2N - 1 + \sum_{i=1}^{N-1} r_i}  \leq  C \sigma_{r+1}.
\end{equation*}
for any cut-off indices $\{r_i\}_{i=1}^{N-1}$. 
\end{theorem}
\begin{proof} 
The result follows from Equations \eqref{eq:disc_sub_int_deco}, \eqref{eq:disc_approx} and Theorem \ref{thm:svd_error} by choosing the subspace
\begin{equation*}
X = \mathop{span} \{ w_{h,\Omega_1},\ldots,w_{h,\Omega_N}, w_{h,f1},\ldots,w_{h,f(N-1)}, \tilde{u}_{h,s}(y) \}.
\end{equation*}
Then the dimension of this space is $\dim(X) = N + (N-1) + \sum_{i=1}^{N-1} r_i$.
\end{proof}
In numerical examples, we have observed that the singular values decay exponentially. This decay is not investigated analytically.

\section{Numerical examples} 

\begin{figure}
     \centering
     \begin{subfigure}[t]{0.45\textwidth}
         \centering
         \includegraphics[width=\textwidth]{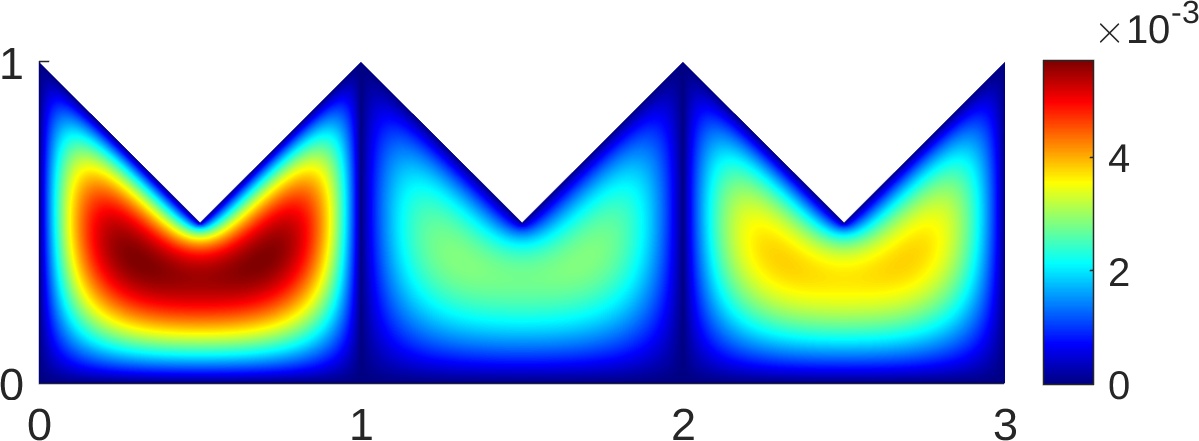}
         \subcaption{The subdomain solutions $u_{h,\Omega_1}(y)$, $u_{h,\Omega_2}(y)$ and $u_{h,\Omega_3}(y)$ in the same picture.}
         \label{fig:domain_sols}
     \end{subfigure}
     \hfill
     \begin{subfigure}[t]{0.45\textwidth}
         \centering
         \includegraphics[width=\textwidth]{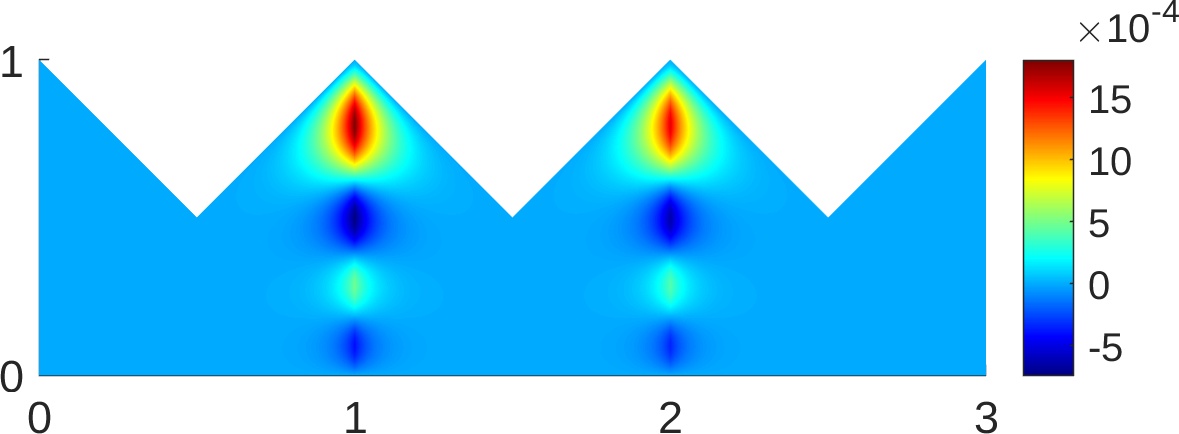}
         \subcaption{The fast solutions $u_{h,f1}(y)$ and $u_{h,f2}(y)$  for the subdomain interfaces.}
         \label{fig:fast_sols}
     \end{subfigure}
     \hfill\\
     \begin{subfigure}[t]{0.45\textwidth}
         \centering
         \includegraphics[width=\textwidth]{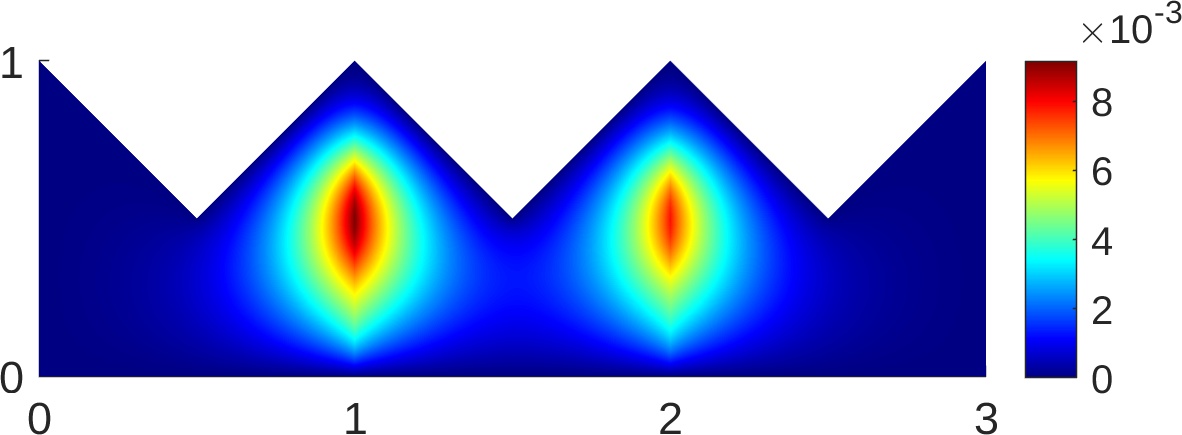}
         \subcaption{The slow solution $\tilde{u}_{h,s}(y)$ for the subdomain interfaces.}
         \label{fig:slow_sols}
     \end{subfigure}
     \hfill
     \begin{subfigure}[t]{0.45\textwidth}
         \centering
         \includegraphics[width=\textwidth]{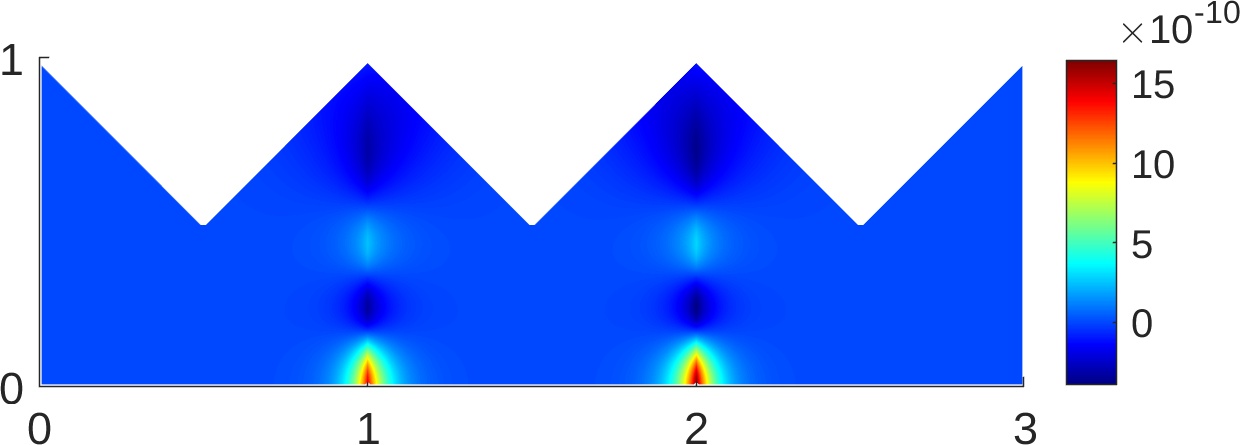}
         \subcaption{The relative error between our approximation and the FE solution calculated for the specific $y$.}
         \label{fig:error}
     \end{subfigure}
     \hfill
     \caption{The different components of the solution and the error for the geometry of Fig. \ref{fig:geom_crown} with $f=x_2$ and $y=[2,4,3]$. The sum of the intermediate solutions \ref{fig:domain_sols}, \ref{fig:fast_sols} and \ref{fig:slow_sols} approximates the true FE solution to a very high accuracy.}
     \label{fig:3domains}
\end{figure}

We numerically implemented the method described in Section \ref{sec:fem}. The examples below use the same geometry as in Fig. \ref{fig:geom_crown} with $N=3$, $f=x_2$ and linear finite element space with dimension $18,785$. However, the layer shape, as well as the value of $N$ and the function $f$ could all easily be modified. First, the different components of the solution for the parameter vector $y=[2,4,3]$ are visualised. The subdomain solution components $u_{h,\Omega_i}(y)=y_i^{-1}w_{h,\Omega_i}$ for $i=1,2,3$ are given in Fig. \ref{fig:domain_sols}, the fast components of the interface solutions $u_{h,fi}(y)=\frac{2}{y_i+y_{i+1}}w_{h,fi}$ for $i=1,2$ are in Fig. \ref{fig:fast_sols} and the slow component $\tilde{u}_{h,s}(y)$ is in Fig. \ref{fig:slow_sols}.

The cut-off indices $\{r_i\}_{i=1}^{N-1}$ that define the discrete slow-fast decomposition were chosen in such a way that $\sigma^{(i)}_{r_i} > 10^{-7}$ and $\sigma^{(i)}_{r_i+1} < 10^{-7}$ for $i=1,\ldots,N-1$. The resulting slow subspace was 6-dimensional, i.e. it had three components per one subdomain interface. In other words, $r_i=r=3$ for $i=\{1,2\}$ and the total dimension of the parameter-to-solution map was 11. The sum of the three intermediate solutions approximates the exact solution. We also plotted the relative error between our approximation and the FE solution for the specific parameter $y$. From Fig. \ref{fig:error} we see that the error is nowhere larger than $2\cdot10^{-9}$. The error forms wave-like patterns on interfaces $\Gamma_i$ and decays quickly outside the interface.

\begin{figure}
    \centering
    \includegraphics[width=0.5\textwidth]{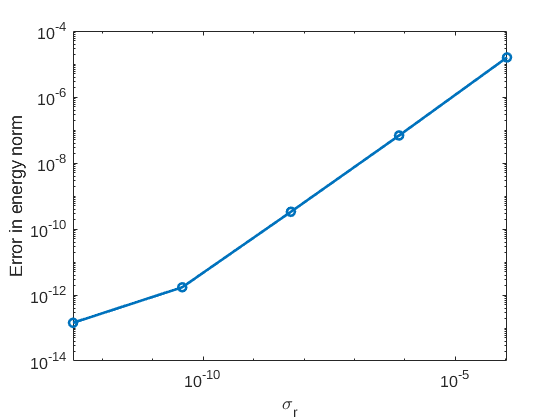}
    \caption{The relative errors in energy norm as a function of the singular value cut-off index $\sigma_r$ for the geometry of Fig. \ref{fig:geom_crown} with $f=x_2$ and $N=3$. We tested multiple parameter vectors $y=[y_1,y_2,y_3]$ where $y_i$ were uniformly sampled from $[1,10]$. The maximum obtained errors for each $\sigma_r$ are plotted.}
    \label{fig:err_plot}
\end{figure}

Relative errors in energy norm are plotted in Fig. \ref{fig:err_plot} as a function of $\sigma_r$. There, we plotted the largest errors obtained for any parameter vector $y=[y_1,y_2,y_3]$ from our sample set. The sample set of parameter vectors was obtained by uniformly sampling $y_i$ from the interval $[1,10]$. The data points correspond to the values $r=[5,4,3,2,1]$ going from left to right. The total dimension of the parameter-to-solution map for each approximation is then $2N + (N-1) r-1$ where $N=3$. We see that the singular values decay exponentially and that the errors are always roughly an order of magnitude smaller than $\sigma_r$. We also noticed that the errors behaved very similarly for different parameters $y$. At slow subspace dimension corresponding to $r=5$ the accuracy compared to the full FE solution is already close to numerical accuracy.

\section{Conclusions and future work}

In this work the parameter dependence of a parametric diffusion equation with a piecewise constant parameter consisting of $N$ layered components was analysed. It was shown that in the case $N=2$ the solution splits into three independent components with an explicit dependence on the parameter $y$. This observation can then be extended for any $N>2$, giving a very efficient way to approximate the solutions for the equation. In the continuous case, we prove error estimates for rectangular domains and in the finite element case for arbitrary mirror symmetric domains. The results also gave upper bounds for the related Kolmogorov $n$-widths.

The authors believe that these techniques could be extended to other geometries, such as $N\times M$ parameter grids. On top of the subdomain-interface decomposition presented in this paper, this would require the treatment of intersection points of four subdomains. We hope that this work sheds light into the reasons behind the effectiveness of sampling-based reduced basis methods for the diffusion equation. The authors also wish to develop reduced basis methods explicitly taking advantage of the solution manifold structure demonstrated in this work.

\bibliographystyle{plain}
\bibliography{refs}

\end{document}